# The topology of deformation spaces of Kleinian groups

By James W. Anderson, Richard D. Canary, and Darryl McCullough*

## Abstract

Let $M$ be a compact, hyperbolizable 3-manifold with nonempty incompressible boundary and let $AH(\pi_1(M))$ denote the space of (conjugacy classes of) discrete faithful representations of $\pi_1(M)$ into $\mathrm{PSL}_2(\mathbf{C})$. The components of the interior $MP(\pi_1(M))$ of $AH(\pi_1(M))$ (as a subset of the appropriate representation variety) are enumerated by the space $\mathcal{A}(M)$ of marked homeomorphism types of oriented, compact, irreducible 3-manifolds homotopy equivalent to $M$. In this paper, we give a topological enumeration of the components of the closure of $MP(\pi_1(M))$ and hence a conjectural topological enumeration of the components of $AH(\pi_1(M))$. We do so by characterizing exactly which changes of marked homeomorphism type can occur in the algebraic limit of a sequence of isomorphic freely indecomposable Kleinian groups. We use this enumeration to exhibit manifolds $M$ for which $AH(\pi_1(M))$ has infinitely many components.

## 1. Introduction

In this paper, we begin a study of the global topology of deformation spaces of Kleinian groups. The basic object of study is the space $AH(\pi_1(M))$ of marked hyperbolic 3-manifolds homotopy equivalent to a fixed compact 3-manifold $M$. The interior $MP(\pi_1(M))$ of $AH(\pi_1(M))$ is very well understood due to work of Ahlfors, Bers, Kra, Marden, Maskit, Sullivan and Thurston. In particular, the components of $MP(\pi_1(M))$ are enumerated by topological data, namely the set $\mathcal{A}(M)$ of marked, compact, oriented, irreducible 3-manifolds homotopy equivalent to $M$, while each component is parametrized by analytic data coming from the conformal boundaries of the hyperbolic 3-manifolds.

Thurston's Ending Lamination Conjecture provides a conjectural classification for elements of $AH(\pi_1(M))$ by data which are partially topologi-

*The research of the second author was supported in part by the National Science Foundation and a fellowship from the Sloan Foundation.

cal, specifically the marked homeomorphism type of the marked hyperbolic 3-manifold $N$ as an element of $\mathcal{A}(M)$, and partially geometric, coming from the conformal boundary of $N$ and the geodesic laminations which encode the asymptotic geometry of any geometrically infinite ends of $N$. However, the data in this conjectural classification do not vary continuously so they do not provide a clear conjectural picture of the global topology. Moreover, the Ending Lamination Conjecture is only known to hold for limited classes of Kleinian groups (see Minsky [31], [32]). In fact, surprisingly little is known about the global topology of $AH(\pi_1(M))$.

In the case that $\pi_1(M)$ is freely indecomposable, the present investigation allows us to give an enumeration of the components of the closure $\overline{MP(\pi_1(M))}$ of $MP(\pi_1(M))$. Since it is conjectured that $\overline{MP(\pi_1(M))} = AH(\pi_1(M))$, this gives a conjectural enumeration of the components of $AH(\pi_1(M))$. In particular, we characterize exactly which components of $MP(\pi_1(M))$ have intersecting closures, by analyzing exactly which changes in marked homeomorphism type can occur in the algebraic limit of a sequence of homotopy equivalent marked hyperbolic 3-manifolds.

One can also think of our work as a study of how the topological data in the Ending Lamination Conjecture vary over $AH(\pi_1(M))$. It follows from our earlier work [4] that these topological data, the marked homeomorphism type, *need not be locally constant* in general. In this paper, we show that marked homeomorphism type *is* locally constant modulo primitive shuffles. Roughly, primitive shuffles are homotopy equivalences obtained by "shuffling" or "rearranging" the way in which the manifold is glued together along the solid torus components of its characteristic submanifold. We hope that further study of how the data in Thurston's Ending Lamination Conjecture vary over $AH(\pi_1(M))$ will eventually yield a conjectural picture of the global topology of $AH(\pi_1(M))$.

In the Sullivan dictionary (see [39]) between Kleinian groups and rational maps, there is an analogy between our work and the study of the parametrization of the Mandelbrot set, or more generally with the study of deformation spaces of rational maps of a fixed type. In this dictionary, the components of $MP(\pi_1(M))$ play the role of hyperbolic components of the Mandelbrot set. Again one may combinatorially enumerate the hyperbolic components of the Mandelbrot set; each component is parametrized by analytic data, in this case the multiplier of the attracting fixed point. The Landing Theorem of Douady and Hubbard may be used to give a complete understanding of which hyperbolic components of the Mandelbrot set have intersecting closures (see Milnor [30] or Schleicher [38]).

It is interesting to contrast the behavior in the two situations. In the Mandelbrot set, one component may "bump" infinitely many other components in the sense that their closures intersect, while any component of $MP(\pi_1(M))$ may bump only finitely many other components. Moreover, the hyperbolic components can accumulate at points in the Mandelbrot set, while the components of $MP(\pi_1(M))$ cannot accumulate in $AH(\pi_1(M))$. On the other hand, in the Mandelbrot set, any point is in the closure of at most two hyperbolic components and any two hyperbolic components bump at at most one point, while arbitrarily many components of $MP(\pi_1(M))$ may bump at a single point and the intersection of the closures of two components of $MP(\pi_1(M))$ can be quite large (see Holt [16]). The global theory of the topology of deformation spaces of hyperbolic 3-manifolds is still in its infancy, but we hope that its study will prove as rich and as rewarding as has the study of the Mandelbrot set.

## 2. Statement of results

In order to state our results precisely, we must introduce some terminology. We will say that a compact, oriented 3-manifold $M$ is *hyperbolizable* if there exists a hyperbolic 3-manifold homeomorphic to the interior of $M$. For a compact, hyperbolizable 3-manifold $M$, let $\mathcal{D}(\pi_1(M))$ be the space of all discrete, faithful representations of $\pi_1(M)$ into $\mathrm{PSL}_2(\mathbf{C})$. Let $\mathrm{Hom}_T(\pi_1(M), \mathrm{PSL}_2(\mathbf{C}))$ denote the set of representations $\rho \in \mathrm{Hom}(\pi_1(M), \mathrm{PSL}_2(\mathbf{C}))$ with the property that $\rho(g)$ is parabolic if $g$ lies in a rank two free abelian subgroup of $\pi_1(M)$. $\mathcal{D}(\pi_1(M))$ is a closed subset of $\mathrm{Hom}_T(\pi_1(M), \mathrm{PSL}_2(\mathbf{C}))$ (see Jørgenson [20]). Set $AH(\pi_1(M)) = \mathcal{D}(\pi_1(M))/\mathrm{PSL}_2(\mathbf{C})$ (where $\mathrm{PSL}_2(\mathbf{C})$ acts by conjugation) and let
$$X_T(\pi_1(M)) = \mathrm{Hom}_T(\pi_1(M), \mathrm{PSL}_2(\mathbf{C}))//\mathrm{PSL}_2(\mathbf{C}),$$
denote the Mumford quotient of $\mathrm{Hom}_T(\pi_1(M), \mathrm{PSL}_2(\mathbf{C}))$ by $\mathrm{PSL}_2(\mathbf{C})$. The space $X_T(\pi_1(M))$ is an algebraic variety, sometimes called the character variety, and $AH(\pi_1(M))$ embeds in $X_T(\pi_1(M))$. (See Chapter V of Morgan-Shalen [34] or Chapter 4 of Kapovich [22] for more details on the character variety.)

Let $MP(\pi_1(M))$ be the subset of $AH(\pi_1(M))$ consisting of minimally parabolic representations with geometrically finite image. Recall that an element $\rho$ of $AH(\pi_1(M))$ is *minimally parabolic* provided that $\rho(g) \in \rho(\pi_1(M))$ is parabolic if and only if $g$ lies in a rank two free abelian subgroup of $\pi_1(M)$. Marden's Quasiconformal Stability Theorem (see [24]) asserts that $MP(\pi_1(M))$ is an open subset of $X_T(\pi_1(M))$. Sullivan [40] showed that the interior of $AH(\pi_1(M))$, as a subset of $X_T(\pi_1(M))$, lies in $MP(\pi_1(M))$, so that $MP(\pi_1(M))$ is the interior of $AH(\pi_1(M))$. It is conjectured that $MP(\pi_1(M))$ is dense in $AH(\pi_1(M))$, see Bers [7], Sullivan [40] or Thurston [41].

The topological type of a hyperbolic 3-manifold $N$ is encoded by its compact core. Recall that a compact 3-dimensional submanifold $C$ of $N$ is a *compact core* for $N$ if the inclusion of $C$ into $N$ is a homotopy equivalence.

A result of Scott [37] implies that every hyperbolic 3-manifold with finitely generated fundamental group has a compact core, and by [28] this core is unique up to homeomorphism. Bonahon [8] proved that if $N$ has incompressible boundary, then $N$ is homeomorphic to the interior of $C$. It is conjectured that every hyperbolic 3-manifold with finitely generated fundamental group is homeomorphic to the interior of its compact core.

If $M$ is a compact, oriented, hyperbolizable 3-manifold, let $\mathcal{A}(M)$ denote the set of *marked homeomorphism types of compact, oriented 3-manifolds homotopy equivalent to $M$*. Explicitly, $\mathcal{A}(M)$ is the set of equivalence classes of pairs $(M', h')$, where $M'$ is a compact, oriented, hyperbolizable 3-manifold and $h' \colon M \to M'$ is a homotopy equivalence, and where two pairs $(M_1, h_1)$ and $(M_2, h_2)$ are equivalent if and only if there exists an orientation-preserving homeomorphism $j \colon M_1 \to M_2$ such that $j \circ h_1$ is homotopic to $h_2$. We denote the class of $(M', h')$ in $\mathcal{A}(M)$ by $[(M', h')]$.

We will use elements of $\mathcal{A}(M)$ to encode the marked homeomorphism type of a marked hyperbolic 3-manifold. For $\rho \in AH(\pi_1(M))$, let $M_\rho$ be a compact core for $N_\rho = \mathbf{H}^3/\rho(\pi_1(M))$ and let $h_\rho \colon M \to M_\rho$ be a homotopy equivalence such that $(h_\rho)_* \colon \pi_1(M) \to \pi_1(M_\rho)$ is conjugate to $\rho$. It is an immediate consequence of Theorem 1 in McCullough, Miller, and Swarup [28] that if $(M_\rho, h_\rho)$ and $(M'_\rho, h'_\rho)$ are two pairs constructed as above, then there exists an orientation-preserving homeomorphism $j \colon M_\rho \to M'_\rho$ such that $j \circ h_\rho$ is homotopic to $h'_\rho$. Therefore, the map $\Psi \colon AH(\pi_1(M)) \to \mathcal{A}(M)$ given by $\Psi(\rho) = [(M_\rho, h_\rho)]$ is well-defined.

Marden's Isomorphism Theorem [24] implies that two elements $\rho_1$ and $\rho_2$ of $MP(\pi_1(M))$ lie in the same component of $MP(\pi_1(M))$ if and only if $\Psi(\rho_1) = \Psi(\rho_2)$, and the Geometrization Theorem of Thurston (see Morgan [33] or Otal [36]) implies that the restriction of $\Psi$ to $MP(\pi_1(M))$ is surjective. Hence, the components of $MP(\pi_1(M))$ are in a one-to-one correspondence with elements of $\mathcal{A}(M)$; the reader is directed to [12] for complete details.

In a previous paper [4], we showed that $\Psi$ need not be locally constant on $AH(\pi_1(M))$. We produced examples $M_k$ (one for each $k \geq 3$) which are obtained by gluing a collection of $k$ $I$-bundles to a solid torus along a collection of parallel annuli. In these examples, all the elements of $\mathcal{A}(M_k)$ are obtained from $M_k$ by "rearranging" or "shuffling" the way in which the $I$-bundles are glued to the solid torus. We showed that any two components of $MP(\pi_1(M_k))$ have intersecting closures. In particular, we showed that there is a sequence of homeomorphic marked hyperbolic 3-manifolds, which converge (algebraically) to a marked hyperbolic 3-manifold which is homotopy equivalent, but not homeomorphic to the approximates. In this paper, we show that, in the case that $M$ has incompressible boundary, a generalization of the phenomenon described in [4] is responsible for all changes of homeomorphism type in the algebraic limit.

Suppose that $M$ has nonempty incompressible boundary. In separate works, Jaco and Shalen [18] and Johannson [19] showed that there exists a *characteristic submanifold* $\Sigma(M)$ of $M$, well-defined up to isotopy, which consists of a disjoint collection of $I$-bundles and Seifert-fibered submanifolds; loosely speaking, this characteristic submanifold captures all the essential annuli and tori in $M$. If $M$ is hyperbolizable, each Seifert-fibered component of $\Sigma(M)$ is homeomorphic to either a solid torus or a thickened torus. A solid torus component $V$ of $\Sigma(M)$ is *primitive* if each component of $\partial M \cap V$ is an annulus whose inclusion into $V$ is a homotopy equivalence. The characteristic submanifold is described in more detail in Section 6.

Given two irreducible 3-manifolds $M_1$ and $M_2$ with nonempty incompressible boundary, a homotopy equivalence $h: M_1 \to M_2$ is a *primitive shuffle* if there exists a finite collection $\mathcal{V}_1$ of primitive solid torus components of $\Sigma(M_1)$ and a finite collection $\mathcal{V}_2$ of solid torus components of $\Sigma(M_2)$, so that $h^{-1}(\mathcal{V}_2) = \mathcal{V}_1$ and so that $h$ restricts to an orientation-preserving homeomorphism from the closure of $M_1 - \mathcal{V}_1$ to the closure of $M_2 - \mathcal{V}_2$. (In Section 5, we will discuss a more general notion of shuffling.) Two elements $[(M_1, h_1)]$ and $[(M_2, h_2)]$ of $\mathcal{A}(M)$ are said to be *primitive shuffle equivalent* if there exists a primitive shuffle $\phi: M_1 \to M_2$ such that $[(M_2, h_2)] = [(M_2, \phi \circ h_1)]$. In Section 7 we observe that if $M$ is hyperbolizable, this gives an equivalence relation on $\mathcal{A}(M)$ and that the resulting quotient map $q: \mathcal{A}(M) \to \widehat{\mathcal{A}}(M)$ is finite-to-one.

Our first main result shows that even though $\Psi$ need not be locally constant, $\widehat{\Psi} = q \circ \Psi$ is locally constant; i.e., marked homeomorphism type is locally constant modulo primitive shuffles. We provide an outline of the proof in Section 8.

THEOREM A. *Let $M$ be a compact, hyperbolizable 3-manifold with nonempty incompressible boundary. Suppose that a sequence $\{\rho_i\} \subset AH(\pi_1(M))$ converges to $\rho \in AH(\pi_1(M))$. Then, for all sufficiently large $i$, $\Psi(\rho_i)$ is primitive shuffle equivalent to $\Psi(\rho)$.*

Our second main result asserts that if the marked homeomorphism types associated to two components of $MP(\pi_1(M))$ are primitive shuffle equivalent, then they have intersecting closures. More specifically, we produce a sequence of marked hyperbolic 3-manifolds in the first component of $MP(\pi_1(M))$ (hence all with the same marked homeomorphism type) which converges algebraically to a geometrically finite marked hyperbolic 3-manifold which has the same marked homeomorphism type as elements of the second component (and hence lies in the boundary of the second component).

THEOREM B. *Let $M$ be a compact, hyperbolizable 3-manifold with nonempty incompressible boundary, and let $[(M_1, h_1)]$ and $[(M_2, h_2)]$ be two elements of $\mathcal{A}(M)$. If $[(M_2, h_2)]$ is primitive shuffle equivalent to $[(M_1, h_1)]$, then the associated components of $MP(\pi_1(M))$ have intersecting closures.*

By combining Theorems A and B, we see that two components of $MP(\pi_1(M))$ have intersecting closures if and only if their corresponding marked homeomorphism types differ by a primitive shuffle.

COROLLARY 1. *Let $M$ be a compact, hyperbolizable 3-manifold with nonempty incompressible boundary, and let $[(M_1, h_1)]$ and $[(M_2, h_2)]$ be two elements of $\mathcal{A}(M)$. The associated components of $MP(\pi_1(M))$ have intersecting closures if and only if $[(M_2, h_2)]$ is primitive shuffle equivalent to $[(M_1, h_1)]$.*

Hence, we can enumerate the components of the closure $\overline{MP(\pi_1(M))}$ of $MP(\pi_1(M))$.

COROLLARY 3. *Let $M$ be a compact, hyperbolizable 3-manifold with nonempty incompressible boundary. Then, the components of $\overline{MP(\pi_1(M))}$ are in a one-to-one correspondence with the elements of $\widehat{\mathcal{A}}(M)$.*

One may combine Corollary 3 with the analysis in [12] to determine exactly when $\overline{MP(\pi_1(M))}$ has infinitely many components. Moreover, we can use Theorem A to establish the existence of 3-manifolds $M$ for which $AH(\pi_1(M))$ has infinitely many components. Recall that it is conjectured that $\overline{MP(\pi_1(M))} = AH(\pi_1(M))$.

A compact, hyperbolizable 3-manifold $M$ with nonempty incompressible boundary has *double trouble* if there exists a toroidal component $T$ of $\partial M$ and homotopically nontrivial simple closed curves $C_1$ in $T$ and $C_2$ and $C_3$ in $\partial M - T$ such that $C_2$ and $C_3$ are not homotopic in $\partial M$, but $C_1$, $C_2$ and $C_3$ are homotopic in $M$. Equivalently, $M$ has double trouble if and only if there is a component $V$ of its characteristic submanifold which is homeomorphic to a thickened torus and whose frontier contains at least two annuli. In [12] it is proven that $\mathcal{A}(M)$ has infinitely many elements if and only if $M$ has double trouble. Since the quotient map $q: \mathcal{A}(M) \to \widehat{\mathcal{A}}(M)$ is finite-to-one, it follows immediately that $\mathcal{A}(M)$ is finite if and only if $\widehat{\mathcal{A}}(M)$ is finite.

COROLLARY 4. *Let $M$ be a compact, hyperbolizable 3-manifold with nonempty incompressible boundary. Then, $\overline{MP(\pi_1(M))}$ has infinitely many components if and only if $M$ has double trouble. Moreover, if $M$ has double trouble, then $AH(\pi_1(M))$ has infinitely many components.*

Another immediate consequence of our main results is that quite often $AH(\pi_1(M))$ is not a manifold. McMullen [29] has used our construction to show that $AH(\pi_1(S))$ is not a manifold if $S$ is a closed hyperbolic surface.

COROLLARY 5. *Let $M$ be a compact, hyperbolizable 3-manifold with nonempty incompressible boundary. If $q: \mathcal{A}(M) \to \widehat{\mathcal{A}}(M)$ is not injective, then $AH(\pi_1(M))$ is not a manifold.*

We also note that the components of $MP(\pi_1(M))$ cannot accumulate in $AH(\pi_1(M))$. This result is used in the proof of Corollary 3, but is also of independent interest.

COROLLARY 2. *Let $M$ be a compact, hyperbolizable 3-manifold with nonempty incompressible boundary. Then, the components of $MP(\pi_1(M))$ cannot accumulate in $AH(\pi_1(M))$. In particular, the closure $\overline{MP(\pi_1(M))}$ of $MP(\pi_1(M))$ is the union of the closures of the components of $MP(\pi_1(M))$.*

## 3. Preliminaries

The purpose of this section is to present some of the background material used in the paper. A *Kleinian group* is a discrete subgroup $\Gamma$ of $\mathrm{PSL}_2(\mathbf{C})$, which we view as acting either on hyperbolic 3-space $\mathbf{H}^3$ via isometries or on the Riemann sphere $\overline{\mathbf{C}}$ via Möbius transformations. The action of $\Gamma$ partitions $\overline{\mathbf{C}}$ into the *domain of discontinuity* $\Omega(\Gamma)$, which is the largest open subset of $\overline{\mathbf{C}}$ on which $\Gamma$ acts properly discontinuously, and the *limit set* $\Lambda(\Gamma)$. A Kleinian group is *nonelementary* if its limit set contains at least three points, and is *elementary* otherwise. A Kleinian group is elementary if and only if it is virtually abelian; recall that a group is *virtually abelian* if it contains a finite index abelian subgroup. We refer the reader to Maskit [26] for a more detailed discussion of the theory of Kleinian groups.

3.1. *Types of Kleinian groups.* There are several specific classes of Kleinian groups which play an important role in this paper.

Given a set $X \subset \mathbf{H}^3 \cup \overline{\mathbf{C}}$ and a Kleinian group $\Gamma$, set

$$\mathrm{st}_\Gamma(X) = \{\gamma \in \Gamma \mid \gamma(X) = X\}.$$

If $\Delta$ is a component of the domain of discontinuity $\Omega(\Gamma)$ of a Kleinian group $\Gamma$, then $\mathrm{st}_\Gamma(\Delta)$ is called a *component subgroup* of $\Gamma$. In the case that $\Gamma$ is finitely generated, it follows from the Ahlfors finiteness theorem that $\mathrm{st}_\Gamma(\Delta)$ is finitely generated and that $\Lambda(\mathrm{st}_\Gamma(\Delta)) = \partial \Delta$; see Lemma 2 of Ahlfors [2].

A finitely generated Kleinian group $\Gamma$ whose domain of discontinuity contains exactly two components $\Delta$ and $\Delta'$ is *quasifuchsian* if $\mathrm{st}_\Gamma(\Delta) = \Gamma$ and is *extended quasifuchsian* otherwise. A quasifuchsian group is quasiconformally conjugate to a Fuchsian group by Theorem 2 of Maskit [25], while an extended quasifuchsian group contains a canonical quasifuchsian subgroup of index two, namely the subgroup stabilizing each of the components of its domain of discontinuity. In particular, the limit set $\Lambda(\Gamma)$ of a quasifuchsian or extended quasifuchsian group is a Jordan curve.

A *web* group is a finitely generated Kleinian group whose domain of discontinuity has at least three components and each of whose component subgroups is quasifuchsian. In particular, each component of the domain of discontinuity of a web group is bounded by a Jordan curve.

A *generalized web* group is either a quasifuchsian, extended quasifuchsian, or web group. The generalized web groups are precisely the finitely generated Kleinian groups with nonempty domain of discontinuity for which every component of the domain of discontinuity is a Jordan domain.

A finitely generated Kleinian group $\Gamma$ is *degenerate* if both its limit set and its domain of discontinuity are nonempty and connected. In this case, there is a conformal homeomorphism $f\colon \mathbf{H}^2 \to \Omega(\Gamma)$. A degenerate group $\Gamma$ is *without accidental parabolics* if an element of the Fuchsian group $f^{-1}\Gamma f$ is parabolic if and only if the corresponding element of $\Gamma$ is parabolic.

A *precisely invariant system of horoballs* $\mathcal{H}$ for a Kleinian group $\Gamma$ is a $\Gamma$-invariant collection of disjoint open horoballs in $\mathbf{H}^3$ such that each horoball in $\mathcal{H}$ is based at a parabolic fixed point of $\Gamma$ and such that there is a horoball in $\mathcal{H}$ based at every parabolic fixed point of $\Gamma$. It is a consequence of the Margulis Lemma, see Benedetti and Petronio [6] or Maskit [26], that every Kleinian group has a precisely invariant system of horoballs. Set $N = \mathbf{H}^3/\Gamma$ and $N^0 = (\mathbf{H}^3 - \mathcal{H})/\Gamma$. If $\Gamma$ is torsion-free, $N$ is a hyperbolic 3-manifold and each component of $\partial N^0$ is either a torus or an open annulus.

The *convex core* $\mathrm{C}(N)$ of $N = \mathbf{H}^3/\Gamma$ is the quotient, by $\Gamma$, of the convex hull $\mathrm{CH}(\Lambda(\Gamma))$ of its limit set. It is the smallest closed, convex submanifold of $N$ whose inclusion into $N$ is a homotopy equivalence. A hyperbolic 3-manifold with finitely generated fundamental group is *geometrically finite* if its convex core has finite volume. A torsion-free Kleinian group is *geometrically finite* if its quotient manifold is geometrically finite. An end $E$ of $N^0$ is *geometrically finite* if it has a neighborhood $U$ such that $U \cap C(N) = \emptyset$, and is *geometrically infinite* otherwise. Note that a finitely generated Kleinian group $\Gamma$ is geometrically finite if and only if every end of $N^0$ is geometrically finite.

An orientable, irreducible 3-manifold $M$ is *topologically tame* if it is homeomorphic to the interior of a compact 3-manifold. A torsion-free Kleinian group is *topologically tame* if its corresponding 3-manifold $N = \mathbf{H}^3/\Gamma$ is topologically tame.

3.2. *Relative compact cores.* We make extensive use of the *relative compact core* of a hyperbolic 3-manifold, which is a compact core which also keeps track of the topology of the parabolic locus. Given a precisely invariant system $\mathcal{H}$ of horoballs for a torsion-free finitely generated Kleinian group $\Gamma$, a compact submanifold $M$ of $N^0$ is a *relative compact core* if it is a compact core for $N^0$ and the intersection of $M$ with each component $Q$ of $\partial N^0$ is a compact core for $Q$. In particular, each component of $P = M \cap \partial N^0$ is either a torus component of $\partial N^0$ or a compact incompressible annulus in an annular component of $\partial N^0$. The fact that every hyperbolic 3-manifold with finitely generated fundamental group possesses a relative compact core follows from work of McCullough [27], see also Kulkarni and Shalen [23]. We often speak of the manifold pair $(M, P)$ as a *relative compact core* of $N$ or of $N^0$.

The relative compact core of a hyperbolic 3-manifold is always a *pared 3-manifold*, several properties of which are described in Lemma 3.1; see [12] or Morgan [33] for a detailed discussion of pared 3-manifolds.

LEMMA 3.1. *Let $N$ be a hyperbolic 3-manifold with finitely generated fundamental group and let $(M, P)$ be a relative compact core for $N$; then, the following hold*:

1. *If $P_1$ is a component of $P$, then $\pi_1(P_1)$ injects into $\pi_1(M)$ and $\pi_1(P_1)$ is a maximal abelian subgroup of $\pi_1(M)$;*

2. *Every noncyclic abelian subgroup of $\pi_1(M)$ is conjugate into the fundamental group of a component of $P$; and*

3. *Every map $\phi\colon (S^1 \times I, S^1 \times \partial I) \to (M, P)$ that induces an injection from $\pi_1(S^1)$ to $\pi_1(M)$ is homotopic, as a map of pairs, to a map $\phi'\colon (S^1 \times I, S^1 \times \partial I) \to (M, P)$ such that $\phi'(S^1 \times I) \subset P$.*

The relative compact core also encodes geometric information about the hyperbolic 3-manifold. For example, we make use of the following standard criterion which guarantees that a Kleinian group is either a generalized web group or a degenerate group without accidental parabolics. We remark that the converse to Lemma 3.2 holds if $\Gamma$ is geometrically finite.

LEMMA 3.2. *Let $\Gamma$ be a finitely generated, torsion-free Kleinian group whose domain of discontinuity $\Omega(\Gamma)$ is nonempty, and let $(M, P)$ be a relative compact core for $N = \mathbf{H}^3/\Gamma$. If every component of $\partial M - P$ is incompressible and if there is no essential annulus in $M$ with one boundary component in $\partial M - P$ and the other in $P$, then $\Gamma$ is either a generalized web group or a degenerate group without accidental parabolics.*

Here, an essential annulus is a properly embedded incompressible annulus which is not properly homotopic into $\partial M$ or, equivalently, into $P$. The proof of Lemma 3.2 is contained in the discussion in Sections 5 and 6 of Abikoff and Maskit [1], particularly Section 6.5. The result there is stated in the language of accidental parabolics, but the equivalence of their statement and ours follows immediately since an accidental parabolic gives rise to an essential annulus in the relative compact core of the form described in the statement of Lemma 3.2.

We will often restrict ourselves to the situation where $(M, P)$ is a relative compact core of a hyperbolic 3-manifold $N$ and every component of $\partial M - P$ is incompressible. In this case, Bonahon [8] showed that $N$ is topologically tame. It will also be useful to notice that in this case the relative compact core is well-defined up to isotopy (which is a relative version of a result of McCullough, Miller and Swarup [28]).

LEMMA 3.3. *Let $\Gamma$ be a finitely generated, torsion-free Kleinian group, let $\mathcal{H}$ be a precisely invariant system of horoballs for $\Gamma$, and let $N^0 = (\mathbf{H}^3 - \mathcal{H})/\Gamma$. If $(M_1, P_1)$ and $(M_2, P_2)$ are two relative compact cores for $N^0$ and if every component of $\partial M_1 - P_1$ is incompressible, then there exists an ambient isotopy of $N^0$ moving $(M_1, P_1)$ onto $(M_2, P_2)$.*

*Proof of Lemma* 3.3. Under our assumptions, Bonahon's theorem [8] guarantees that every end of $N^0$ has a neighborhood which is a product of a compact surface with the half-line. Hence, there exists a relative compact core $(M_3, P_3)$ for $N^0$ which contains both $M_1$ and $M_2$ in its interior. Since $(M_1, P_1)$ and $(M_2, P_2)$ are both homotopy equivalent, as pairs, to $(N^0, \partial N^0)$, they are homotopy equivalent to one another. Thus, since every component of $\partial M_1 - P_1$ is incompressible, every component of $\partial M_2 - P_2$ is also incompressible; see Proposition 1.2 of Bonahon [8] or Section 5.2 of [12].

Now let $i$ stand for either 1 or 2. Since $P_i$ and $P_3$ are compact cores for $\partial N^0$, each component of $\overline{P_3 - P_i}$ is an annulus with one boundary circle lying in $P_i$. Since $M_i$ and $M_3$ are relative compact cores for $N^0$, $M_3$ contains $M_i$ and the frontier of $M_i$ is incompressible, each component $F$ of the frontier of $M_i$ has the property that $\pi_1(F)$ surjects onto $\pi_1(W)$, where $W$ is the component of $\overline{M_3 - M_i}$ which is bounded by $F$. Theorem 10.5 of Hempel [15] then implies that $W$ is a product $F \times I$ with $F = F \times \{0\}$. Since $\overline{P_3 - P_i}$ consists of annuli meeting $P_i$, we may choose the product structure so that $(F \times I) \cap \partial N^0 = \partial F \times I$. It follows that there is an ambient isotopy moving $M_i$ onto $M_3$. Since this is true for both $i = 1$ and $i = 2$, $M_1$ is ambient isotopic to $M_2$. □

3.3. *Geometric convergence.* In contrast to the topology of algebraic convergence described in Section 2, there is a second notion of convergence for Kleinian groups which is more closely allied to the geometry of the quotient hyperbolic 3-manifolds. Say that a sequence $\{\Gamma_i\}$ of Kleinian groups converges *geometrically* to $\widehat{\Gamma}$ if every element $\gamma \in \widehat{\Gamma}$ is the limit of a sequence $\{\gamma_i \in \Gamma_i\}$ and if every accumulation point of every sequence $\{\gamma_i \in \Gamma_i\}$ lies in $\widehat{\Gamma}$.

Jørgenson and Marden [21] observed that if $M$ is a compact, hyperbolizable 3-manifold with nonabelian fundamental group and $\{\rho_i\}$ is a sequence in $\mathcal{D}(\pi_1(M))$ converging to $\rho$, then there exists a subsequence of $\{\rho_i\}$, again called $\{\rho_i\}$, so that $\{\rho_i(\pi_1(M))\}$ converges geometrically to a torsion-free Kleinian group $\widehat{\Gamma}$ with $\rho(\pi_1(M)) \subset \widehat{\Gamma}$. A sequence $\{\rho_i\}$ in $\mathcal{D}(\pi_1(M))$ converges *strongly* to $\rho$ if $\{\rho_i\}$ converges to $\rho$ and if $\{\rho_i(\pi_1(M))\}$ converges geometrically to $\rho(\pi_1(M))$. We also make use of the next result, which follows immediately from Corollary 3.9 of Jørgensen and Marden [21].

LEMMA 3.4. *Suppose that $M$ is a compact, hyperbolizable 3-manifold with non-abelian fundamental group, that $\{\rho_i\}$ is a sequence in $\mathcal{D}(\pi_1(M))$ converging to $\rho$, and that $\{\rho_i(\pi_1(M))\}$ converges geometrically to a Kleinian group $\widehat{\Gamma}$. If $\gamma \in \widehat{\Gamma}$ and if there exists $n > 0$ such that $\gamma^n \in \rho(\pi_1(M))$, then $\gamma \in \rho(\pi_1(M))$.*

The following lemma highlights the geometric significance, on the level of the quotient manifolds, of the geometric convergence of a sequence of Kleinian groups. Setting notation, let 0 denote a choice of basepoint for $\mathbf{H}^3$, and let $p_i\colon \mathbf{H}^3 \to N_i = \mathbf{H}^3/\Gamma_i$ and $p\colon \mathbf{H}^3 \to \widehat{N} = \mathbf{H}^3/\widehat{\Gamma}$ be the covering maps. Let $B_R(0) \subset \mathbf{H}^3$ be a ball of radius $R$ centered at the basepoint 0.

LEMMA 3.5 (Theorem E.1.13 in Benedetti-Petronio [6]). *A sequence of torsion-free Kleinian groups $\{\Gamma_i\}$ converges geometrically to a torsion-free Kleinian group $\widehat{\Gamma}$ if and only if there exists a sequence $\{(R_i, K_i)\}$ and a sequence of orientation-preserving maps $\widetilde{f}_i\colon B_{R_i}(0) \to \mathbf{H}^3$ such that the following hold:*

1. *$R_i \to \infty$ and $K_i \to 1$ as $i \to \infty$;*
2. *The map $\widetilde{f}_i$ is a $K_i$-bilipschitz diffeomorphism onto its image, $\widetilde{f}_i(0) = 0$, and $\{\widetilde{f}_i|_A\}$ converges, in the $C^\infty$-topology, to the identity on any compact subset $A$ of $\mathbf{H}^3$; and*
3. *$\widetilde{f}_i$ descends to a map $f_i\colon Z_i \to N_i$, where $Z_i = B_{R_i}(0)/\widehat{\Gamma}$ is a submanifold of $\widehat{N}$; moreover, $f_i$ is also an orientation-preserving $K_i$-bilipschitz diffeomorphism onto its image.*

*Remark.* In the actual statement of Theorem E.1.13 in [6] it is only required in (2) that $\{\widetilde{f}_i|_A\}$ converges, in the $C^\infty$-topology, to the identity on any compact subset $A$ of $\mathbf{H}^3$. Also, it is only required in (3) that each $\widetilde{f}_i$ descends to a map $f_i\colon Z_i \to N_i$. In remarks E.1.12 and E.1.19 it is noted that our stronger formulations of conclusions (2) and (3) follow from their statement.

The following lemma, whose proof is the same as that of Proposition 3.3 of Canary and Minsky [13] (see also Lemma 7.2 of [3]), allows us, in the case of a geometric limit of an algebraically convergent sequence of Kleinian groups, to see the relationship between the diffeomorphisms $\{f_i\}$ and the representations $\{\rho_i\}$.

LEMMA 3.6. *Suppose that $M$ is a compact, hyperbolizable 3-manifold, that $\{\rho_i\}$ is a sequence in $\mathcal{D}(\pi_1(M))$ converging to $\rho$, and that $\{\rho_i(\pi_1(M))\}$ converges geometrically to a Kleinian group $\widehat{\Gamma}$. Let $N = \mathbf{H}^3/\rho(\pi_1(M))$ and $\widehat{N} = \mathbf{H}^3/\widehat{\Gamma}$, let $\pi\colon N \to \widehat{N}$ be the covering map, and let $\{f_i\colon Z_i \to N_i\}$ be* *the sequence of bilipschitz diffeomorphisms produced by Lemma* 3.5. *Suppose that $K$ is a compact subset of $N$ such that $\pi_1(K)$ injects into $\pi_1(N)$ and let $\pi': K \to \pi(K)$ denote the restriction of $\pi$ to $K$. Then, for all sufficiently large $i$, $(f_i \circ \pi')_*$ agrees with the restriction of $\rho_i \circ \rho^{-1}$ to $\pi_1(K)$, where both are regarded as giving maps of $\pi_1(K) \subset \pi_1(N) = \rho(\pi_1(M))$ into $\pi_1(N_i) = \rho_i(\pi_1(M))$.*

We also need to make use of two results from [3]. The first, Proposition 4.2, shows that any "new" element of the geometric limit must map the limit set of a topologically tame subgroup of the algebraic limit off of the the limit set of any other topologically tame subgroup of the algebraic limit, except perhaps for one parabolic fixed point which they may have in common. Given a pair $\Phi$ and $\Phi'$ of subgroups of a Kleinian group $\Gamma$, let $P(\Phi, \Phi')$ be the set of points $x \in \Lambda(\Gamma)$ such that $\text{st}_\Phi(x)$ and $\text{st}_{\Phi'}(x)$ are both rank one parabolic subgroups and $\langle \text{st}_\Phi(x), \text{st}_{\Phi'}(x) \rangle$ is a rank two parabolic subgroup of $\Gamma$.

PROPOSITION 3.7. *Suppose that $M$ is a compact hyperbolizable $3$-manifold, that $\{\rho_i\}$ is a sequence in $\mathcal{D}(\pi_1(M))$ converging to $\rho$, and that $\{\rho_i(\pi_1(M))\}$ converges geometrically to a Kleinian group $\widehat{\Gamma}$. If $\Gamma_1$ and $\Gamma_2$ are two (possibly equal) topologically tame subgroups of $\rho(\pi_1(M))$ and if $\gamma \in \widehat{\Gamma} - \rho(\pi_1(M))$, then $P(\Gamma_1, \gamma\Gamma_2\gamma^{-1}) = \emptyset$. Moreover, $\Lambda(\Gamma_1 \cap \gamma\Gamma_2\gamma^{-1}) = \Lambda(\Gamma_1) \cap \gamma(\Lambda(\Gamma_2))$, and $\Lambda(\Gamma_1 \cap \gamma\Gamma_2\gamma^{-1})$ contains at most one point.*

The second result, Proposition 5.2 from [3], is a consequence of Thurston's covering theorem [42], as generalized by Canary [11].

PROPOSITION 3.8. *Suppose that $M$ is a compact hyperbolizable $3$-manifold, that $\{\rho_i\}$ is a sequence in $\mathcal{D}(\pi_1(M))$ converging to $\rho$, and that $\{\rho_i(\pi_1(M))\}$ converges geometrically to a Kleinian group $\widehat{\Gamma}$. Let $N = \mathbf{H}^3/\rho(\pi_1(M))$ and $\widehat{N} = \mathbf{H}^3/\widehat{\Gamma}$, and let $\pi: N \to \widehat{N}$ be the covering map. Let $\mathcal{H}$ be a precisely invariant system of horoballs for $\rho(\pi_1(M))$. Suppose that $\rho(\pi_1(M))$ is topologically tame and that $E$ is a geometrically infinite end of $N^0$. Then, there is a neighborhood $U$ of $E$ such that the restriction of $\pi$ to $U$ is an injection.*

# 4. Precisely embedded subgroups and relative compact carriers

In this section we develop a geometric criterion, expressed in terms of the limit sets, which guarantees that a collection of generalized web subgroups of the fundamental group of a hyperbolic 3-manifold $N$ is associated to a collection of disjoint submanifolds of $N$. This generalizes Proposition 6.1 in [3] which shows that a single precisely embedded generalized web subgroup is associated to a compact submanifold.

A generalized web subgroup $\Theta$ of a Kleinian group $\Gamma$ is *precisely embedded* if it is maximal in the sense that $\text{st}_\Gamma(\Lambda(\Theta)) = \Theta$ and if, for each $\gamma \in \Gamma - \Theta$, there exists a component of $\Omega(\Theta)$ whose closure contains $\gamma(\Lambda(\Theta))$. More generally, a collection $\{\Gamma_1, \ldots, \Gamma_n\}$ of generalized web subgroups of $\Gamma$ is a *precisely embedded system* if each $\Gamma_j$ is a precisely embedded subgroup of $\Gamma$ and if, for $\gamma \in \Gamma$ and $j \neq k$, there is a component of $\Omega(\Gamma_j)$ whose closure contains $\gamma(\Lambda(\Gamma_k))$.

For the remainder of this section, we adopt the following notation. Let $\Gamma$ be a torsion-free Kleinian group, let $N = \mathbf{H}^3/\Gamma$, let $\mathcal{H}$ be a precisely invariant system of horoballs for $\Gamma$, and let $N^0 = (\mathbf{H}^3 - \mathcal{H})/\Gamma$. For a subgroup $\Gamma_j$ of $\Gamma$, let $\mathcal{H}_j$ denote the collection of those horoballs in $\mathcal{H}$ which are based at fixed points of parabolic elements of $\Gamma_j$, let $N_j = \mathbf{H}^3/\Gamma_j$, and let $N_j^0 = (\mathbf{H}^3 - \mathcal{H}_j)/\Gamma_j$. Let $p_j \colon N_j \to N$, $q_j \colon \mathbf{H}^3 \to N_j$, and $\pi \colon \mathbf{H}^3 \to N$ be the covering maps. If $R_j'$ is a relative compact core for $N_j^0$ and $p_j$ is injective on $R_j'$, then we call the image $R_j = p_j(R_j')$ in $N$ a *relative compact carrier* of the subgroup $\Gamma_j$ of $\Gamma$.

The main result of this section, Proposition 4.2, asserts that a precisely embedded system $\{\Gamma_1, \ldots, \Gamma_n\}$ of nonconjugate generalized web subgroups of a torsion-free Kleinian group $\Gamma$ has a system of disjoint relative compact carriers in $\mathbf{H}^3/\Gamma$. We begin by establishing a partial converse to Proposition 4.2, whose proof will serve as a guide for that of the main result. This partial converse demonstrates that the condition of being a precisely embedded system almost characterizes which collections of generalized web subgroups are associated to disjoint collections of compact submanifolds.

LEMMA 4.1. *Let $\Gamma$ be a torsion-free Kleinian group, and let $\{\Gamma_1, \ldots, \Gamma_n\}$ be a collection of generalized web subgroups of $\Gamma$. Suppose there exists a disjoint collection $\{R_1, \ldots, R_n\}$ of submanifolds of $N = \mathbf{H}^3/\Gamma$ such that, for all $j$, $R_j$ is a relative compact carrier for $\Gamma_j$, and if $R_j$ is an I-bundle, then no component of the closure of $N^0 - R_j$ is a compact twisted I-bundle whose associated $\partial$I-bundle lies in $\partial R_j$. Then, $\{\Gamma_1, \ldots, \Gamma_n\}$ is a precisely embedded system of generalized web subgroups of $\Gamma$.*

In the proof of Lemma 4.1 and its converse, Proposition 4.2, we make extensive use of the notion of a spanning disc for a quasifuchsian or extended quasifuchsian subgroup. Roughly, a spanning disc is the lift to the universal cover of a properly embedded surface representing the subgroup. More precisely, a *spanning disc* for a quasifuchsian or extended quasifuchsian subgroup $\Theta$ of $\Gamma$ is a disc $D$ in $\mathbf{H}^3$ which is precisely invariant under $\Theta$ in $\Gamma$ and which extends to a closed disc $\overline{D}$ in $\mathbf{H}^3 \cup \overline{\mathbf{C}}$ with boundary $\Lambda(\Theta)$. (Recall that a subset $X$ of $\mathbf{H}^3 \cup \overline{\mathbf{C}}$ is *precisely invariant* under a subgroup $\Theta$ of $\Gamma$ if $\text{st}_\Gamma(X) = \Theta$ and if $\gamma(X) \cap X = \emptyset$ for all $\gamma \in \Gamma - \Theta$.) A collection $\{D_1, \ldots, D_q\}$ of disjoint discs in $\mathbf{H}^3$ is a *system of spanning discs* for a precisely embedded system of quasifuchsian and extended quasifuchsian subgroups $\{\Theta_1, \ldots, \Theta_q\}$ of $\Gamma$ if each $D_j$ is a spanning disc for $\Theta_j$ and if $\gamma(D_j)$ and $D_k$ are disjoint for all $\gamma \in \Gamma$ and for all $k \neq j$. A key result used in the proof of Proposition 4.2 is that every precisely embedded system of quasifuchsian and extended quasifuchsian subgroups admits a system of spanning discs (see Lemma 6.3 in [3]).

*Proof of Lemma* 4.1. We begin by constructing an extension $S_j$ of each $R_j$ by appending the portions of the cusps of $N$ which "project onto" $R_j$. Let $s\colon N - N^0 \to \partial N^0$ be the map which takes a point in $N - N^0$ to the closest point in $\partial N^0$ (i.e. the map given by perpendicular projection of each cusp onto its boundary). For each $1 \leq j \leq n$, set $S_j = R_j \cup s^{-1}(R_j \cap N^0)$ and let $\widetilde{S}_j$ denote the component of the preimage of $S_j$ in $\mathbf{H}^3$ which is stabilized by $\Gamma_j$. In particular, each component $T'$ of the closure $\overline{\partial R_j - \partial N^0}$ of $\partial R_j - N^0$ is a subset of a boundary component $T$ of $S_j$ and $T$ is obtained from $T'$ by "appending a cusp" to each component of $\partial T'$. Notice that $\widetilde{S}_j$ is precisely invariant under $\Gamma_j$ in $\Gamma$. Moreover, since $S_j$ and $S_k$ are disjoint subsets of $N$, $\gamma(\widetilde{S}_j)$ and $\widetilde{S}_k$ are disjoint for all $\gamma \in \Gamma$ and all $j \neq k$.

By assumption, $R_j$ lifts to a relative compact core $R'_j$ of $N^0_j$. If $D$ is a component of $\partial \widetilde{S}_j$, then $\widehat{D} = q_j(D) \cap R'_j$ bounds a neighborhood of an end of $N^0_j$. If $\widehat{D}$ bounds a geometrically finite end of $N^0_j$ then $D$ is a spanning disc for a (necessarily quasifuchsian) component subgroup $\mathrm{st}_{\Gamma_j}(\Delta)$ of $\Gamma_j$. On the other hand, if $\widehat{D}$ bounds a geometrically infinite end of $N^0_j$ and $A$ is the component of $\mathbf{H}^3 - \widetilde{S}_j$ bounded by $D$, then $A \subset \mathrm{CH}(L_\Gamma)$, so that $\overline{A} \cap \overline{\mathbf{C}}$ is contained in $\Lambda(\Gamma_j)$.

We begin by showing that if $\gamma \in \Gamma - \Gamma_j$, then $\gamma(\Lambda(\Gamma_j))$ is contained in the closure of a component of $\Omega(\Gamma_j)$. If $\gamma \in \Gamma - \Gamma_j$, then $\widetilde{S}_j$ and $\gamma(\widetilde{S}_j)$ are disjoint. Let $D$ be the component of $\partial \widetilde{S}_j$ which lies between $\gamma(\widetilde{S}_j)$ and $\widetilde{S}_j$. If $\widehat{D} = q_j(D) \cap R'_j$ bounds a neighborhood of a geometrically infinite end of $N^0_j$, let $\Delta'$ be a component of $\gamma(\Omega(\Gamma_j))$, such that the region $H$ in $\mathbf{H}^3$ between $\gamma(\widetilde{S}_j)$ and $\Delta'$ is contained in the component $A$ of $\mathbf{H}^3 - \widetilde{S}_j$ bounded by $D$. Then $\Delta' \subset \overline{A} \cap \overline{\mathbf{C}} \subset \Lambda(\Gamma_j)$, which is impossible, since $\Lambda(\Gamma_j)$ has empty interior. Therefore $D$ must be a spanning disc for a component subgroup $\mathrm{st}_{\Gamma_j}(\Delta)$ of $\Gamma_j$. Since $\gamma(\widetilde{S}_j)$ is contained entirely between $D$ and $\Delta$, it follows that $\gamma(\Lambda(\Gamma_j)) \subset \overline{\Delta}$ as desired.

Now suppose that some $\gamma \in \Gamma - \Gamma_j$ fixes $\Lambda(\Gamma_j)$. We see from the above that there exists a component $\Delta$ of $\Omega(\Gamma_j)$ such that $\gamma(\Lambda(\Gamma_j))$ lies in $\overline{\Delta}$. Thus $\Lambda(\Gamma_j)$ must be equal to $\partial \Delta$ and hence $\Gamma_j$ is either quasifuchsian or extended quasifuchsian and so $R_j$ is an $I$-bundle. Let $\Theta = \mathrm{st}_\Gamma(\Lambda(\Gamma_j))$ and $N_\Theta = \mathbf{H}^3/\Theta$. Let $\mathcal{H}_\Theta$ denote the horoballs in $\mathcal{H}$ which are based at parabolic fixed points of $\Theta$ and $N^0_\Theta = (\mathbf{H}^3 - \mathcal{H}_\Theta)/\Theta$. Then $R_j$ lifts to a relative compact carrier $R_\Theta$ for $\Gamma_j$ in $\Theta$. Notice that since $\Theta$ is quasifuchsian or extended quasifuchsian, $\Gamma_j$ is finitely generated and $\Lambda(\Gamma_j) = \Lambda(\Theta)$, $\Gamma_j$ must have finite index in $\Theta$. If $\Theta$ is quasifuchsian, then $N_\Theta^0$ is an untwisted $\mathbf{R}$-bundle over a compact surface and no proper finite index subgroup admits a relative compact carrier. So, $\Theta$ must be extended quasifuchsian, in which case $N_\Theta^0$ is a twisted $\mathbf{R}$-bundle over a compact surface and the only proper finite index subgroup which admits a relative compact carrier is its index two quasifuchsian component subgroup. Thus, one component $B$ of the closure of $N_\Theta^0 - R_\Theta$ is a compact twisted $I$-bundle whose associated $\partial I$-bundle lies in $\partial R_\Theta$, and hence, since $R_\Theta$ is a lift of $R_j$, one may use Theorem 10.5 of Hempel [15] to check that one component of the closure of $N^0 - R_j$ is a compact twisted $I$-bundle whose associated $\partial I$-bundle lies in $\partial R_j$. Since we have explicitly ruled this situation out we may conclude that $\mathrm{st}_\Gamma(\Lambda(\Gamma_j)) = \Gamma_j$ and thus that $\Gamma_j$ is a precisely embedded subgroup of $\Gamma$.

If $\gamma \in \Gamma$ and $j \ne k$, then $\gamma(\widetilde{S}_j)$ and $\widetilde{S}_k$ are disjoint, and so there exists a component $D$ of $\partial \widetilde{R}_k$ which lies between $\gamma(\widetilde{S}_j)$ and $\widetilde{S}_k$. We may argue exactly as before to show that $D$ is a spanning disc of a component subgroup $\mathrm{st}_{\Gamma_j}(\Delta)$ of $\Gamma_j$ and that $\gamma(\Lambda(\Gamma_k))$ lies in the closure of $\Delta$. This completes the proof that $\{\Gamma_1, \ldots, \Gamma_n\}$ is a precisely embedded system of generalized web subgroups of $\Gamma$. $\square$

The proof of Proposition 4.2 consists of reversing the process described in the proof of the above lemma. It closely resembles the argument given in the proof of Lemma 6.1 in [3]. We first construct an equivariant collection of spanning discs for the component subgroups of the $\Gamma_j$. We then observe that the region bounded by these spanning discs is precisely invariant under $\Gamma_j$ in $\Gamma$. Hence, if we construct a relative compact core for $N_j^0$ which lies in the image of this region in $N_j^0$, it embeds in $N^0$ and provides a relative compact carrier for $\Gamma_j$ in $\Gamma$.

PROPOSITION 4.2. *Let $\Gamma$ be a torsion-free Kleinian group, and let $\{\Gamma_1, \ldots, \Gamma_n\}$ be a precisely embedded system of nonconjugate generalized web subgroups of $\Gamma$. Then, there exists a disjoint collection $\{Y_1, \ldots, Y_n\}$ of relative compact carriers for $\{\Gamma_1, \ldots, \Gamma_n\}$ in $N = \mathbf{H}^3/\Gamma$.*

*Proof of Proposition* 4.2. Let $\{\Theta_1', \ldots, \Theta_q'\}$ denote a maximal collection of nonconjugate quasifuchsian subgroups of $\Gamma$ which arise as component subgroups of groups in the collection $\{\Gamma_1, \ldots, \Gamma_n\}$. Since each $\Gamma_j$ is precisely embedded, Lemma 6.2 of [3] asserts that, for all $m$, either $\Theta_m'$ is precisely embedded in $\Gamma$, in which case we let $\Theta_m = \Theta_m'$, or is an index two subgroup of a precisely embedded extended quasifuchsian subgroup $\Theta_m$ of $\Gamma$. Combining this with the fact that $\{\Gamma_1, \ldots, \Gamma_n\}$ is a precisely embedded system of generalized web subgroups of $\Gamma$, it is easy to check that $\{\Theta_1, \ldots, \Theta_q\}$ is a precisely embedded system of quasifuchsian and extended quasifuchsian subgroups of $\Gamma$. Lemma 6.3 of [3] then implies that there exists a system $\{D_1, \ldots, D_q\}$ of spanning discs for $\{\Theta_1, \ldots, \Theta_q\}$ such that each $F_m = \pi(D_m)$ is a properly embedded surface in $N$. We may further assume that the intersection of each $F_m$ with each component of $N - N^0$ is totally geodesic.

Let $\{B_1, \ldots, B_q\}$ be a collection of pairwise disjoint, closed regular neighborhoods of the surfaces $\{F_1, \ldots, F_q\}$ in $N$. As above, we may assume that each component of the intersection of each $B_m$ with $N - N^0$ is bounded by a pair of totally geodesic surfaces. If $\Theta_m$ is quasifuchsian, then $F_m$ is orientable and $B_m$ is a standard $I$-bundle over $F_m$. If $\Theta_m$ is extended quasifuchsian, then $F_m$ is nonorientable and $B_m$ is a twisted $I$-bundle. Let $\widetilde{B}_m$ be the component of the preimage of $B_m$ in $\mathbf{H}^3$ which contains $D_m$, and note that $\widetilde{B}_m$ is a closed regular neighborhood of $D_m$ which is precisely invariant under $\Theta_m$ in $\Gamma$. Moreover, each component of $\partial \widetilde{B}_m$ is a spanning disc for $\Theta'_m$.

Renumber so that $\Gamma_j$ is quasifuchsian or extended quasifuchsian if $1 \leq j \leq l$, and so that $\Gamma_j$ is a web group if $l+1 \leq j \leq n$. For each quasifuchsian or extended quasifuchsian $\Gamma_j$, we may, perhaps after renumbering $\{\Theta_1, \ldots, \Theta_q\}$, assume that $\Gamma_j = \Theta_j$. In this case, $Y_j = B_j \cap N^0$ is a relative compact carrier for $\Gamma_j$ in $\Gamma$. By construction $Y_j$ is disjoint from $Y_k$ for $j \neq k$ and $1 \leq j, k \leq l$.

Suppose now that $\Gamma_j$ is a web group. For each component $\Delta$ of $\Omega(\Gamma_j)$, there exists a translate $\widetilde{B}_\Delta$ of an element of $\{\widetilde{B}_1, \ldots, \widetilde{B}_q\}$ whose boundary components are spanning discs for $\mathrm{st}_{\Gamma_j}(\Delta)$. Let $\widetilde{X}_\Delta$ denote the component of $\mathbf{H}^3 - \widetilde{B}_\Delta$ which does not contain $\Delta$ in its Euclidean closure in $\mathbf{H}^3 \cup \overline{\mathbf{C}}$. Let $\widetilde{X}_j = \bigcap \widetilde{X}_\Delta$, where the intersection is taken over all components $\Delta$ of $\Omega(\Gamma_j)$. By construction, $\widetilde{X}_j$ is invariant under $\Gamma_j$.

Now suppose that $\widetilde{X}_j$ intersects some translate $\gamma(\widetilde{B}_m)$ of an element of $\{\widetilde{B}_1, \ldots, \widetilde{B}_q\}$. Then some boundary component $\gamma(D)$ of $\gamma(\widetilde{B}_m)$ either intersects the boundary of $\widetilde{X}_j$ transversely or lies in the interior of $\widetilde{X}_j$. The first possibility is ruled out by the fact that all the translates of the elements of $\{\widetilde{B}_1, \ldots, \widetilde{B}_q\}$ are disjoint. The second possibility implies that $\gamma(\Lambda(\Theta'_m)) \subset \Lambda(\Gamma_j)$. However, as $\Theta'_m$ is a component subgroup of some $\Gamma_k$ and $\{\Gamma_1, \ldots, \Gamma_n\}$ is a precisely embedded system of generalized web groups, this implies that $\gamma(\Lambda(\Theta_m))$ bounds a component $\Delta$ of $\Omega(\Gamma_j)$. Since both $\Gamma_j$ and $\Theta_m$ are precisely embedded in $\Gamma$, $\mathrm{st}_{\Gamma_j}(\Delta)$ must be equal to $\gamma\Theta'_m\gamma^{-1}$ and hence $\gamma(\widetilde{B}_m) = \widetilde{B}_\Delta$ is disjoint from $\widetilde{X}_j$ by construction. Therefore, each $\widetilde{X}_j$ must be disjoint from any translate of an element of $\{\widetilde{B}_1, \ldots, \widetilde{B}_q\}$.

Suppose that $\widetilde{X}_j$ is not precisely invariant under $\Gamma_j$ in $\Gamma$. If $\gamma(\widetilde{X}_j) = \widetilde{X}_j$, then $\gamma(\Lambda(\Gamma_j)) = \Lambda(\Gamma_j)$, which implies that $\gamma \in \Gamma_j$ (since $\Gamma_j$ is precisely embedded). Thus, either there exists a boundary component of $\widetilde{X}_j$ which intersects a boundary component of $\gamma(\widetilde{X}_j)$ transversely for some $\gamma \in \Gamma - \Gamma_j$, or there exist a component $D$ of $\partial \widetilde{X}_j$ and an element $\gamma \in \Gamma - \Gamma_j$ so that $\gamma(D)$ lies in the interior of $\widetilde{X}_j$. The former case is ruled out by the construction of $\{\widetilde{B}_1, \ldots, \widetilde{B}_q\}$. In the latter case, $\gamma(D)$ bounds a translate of an element of $\{\widetilde{B}_1, \ldots, \widetilde{B}_q\}$ which must itself intersect $\widetilde{X}_j$, which we ruled out in the above paragraph. Therefore, $\widetilde{X}_j$ is precisely invariant under $\Gamma_j$ in $\Gamma$. A very similar argument gives that $\gamma(\widetilde{X}_j)$ and $\widetilde{X}_k$ are disjoint for all $\gamma \in \Gamma$ and all $j \neq k$ with $l+1 \leq j, k \leq n$.

Suppose that $Q$ is a component of $\partial N_j^0$ and that $H$ is a horoball in $\mathcal{H}$ whose boundary covers $Q$ and which is based at a parabolic fixed point $x \in \Lambda(\Gamma_j)$. If $\Delta$ is a component of $\Omega(\Gamma_j)$ and $x$ is not in $\partial \Delta$, then $\widetilde{B}_\Delta$ cannot intersect $H$; so $H \subset \widetilde{X}_\Delta$. If $x \in \partial \Delta$, then $\widetilde{B}_\Delta$ intersects $H$ in the region between two totally geodesic planes, so that $H \cap \widetilde{X}_\Delta$ is bounded by a totally geodesic hyperplane. Since $x$ lies on the boundary of at most two components of $\Omega(\Gamma_j)$ and the boundaries of the $\widetilde{X}_\Delta$ are disjoint, $H \cap \widetilde{X}_j^0$ must be nonempty. If $x$ is a rank 2 parabolic fixed point, it cannot lie in the boundary of any component of $\Omega(\Gamma_j)$, so that no component of $\widetilde{B}_\Delta$ can intersect $H$, and $Q$ must be contained entirely within $\widetilde{X}_j$. Notice that in all cases, $H \cap \widetilde{X}_j$ is an open topological half-space bounded by a properly embedded topological plane.

Since $\widetilde{X}_j$ is the complement of a disjoint collection of open topological half-spaces in $\mathbf{H}^3$ (which do not accumulate at any point in $\mathbf{H}^3$), it is contractible. Similarly, $\widetilde{X}_j^0 = \widetilde{X}_j - \mathcal{H}$ is contractible, as each component of $\mathcal{H} \cap \widetilde{X}_j$ is an open topological half-space whose boundary is a properly embedded topological plane and the components do not accumulate within $\mathbf{H}^3$. Therefore, the inclusion of $X_j^0 = q_j(\widetilde{X}_j^0)$ into $N_j^0$ is a homotopy equivalence. If $Q$ is an annular component of $\partial N_j^0$, then the argument in the above paragraph shows that $X_j^0 \cap Q$ contains an incompressible annulus, while if $Q$ is a torus component of $\partial N_j^0$ then this shows that $Q \subset X_j^0$. For any component $Q$ of $\partial N_j^0$, let $A_Q$ be a compact core for $\partial Q$ which is contained in $X_j^0$. A theorem of McCullough [27] (see also Kulkarni and Shalen [23]) then implies that we may find a compact core $Y_j'$ for $X_j^0$ so that $Y_j' \cap Q = A_Q$ for every component $Q$ of $\partial N_j^0$. Since the inclusion of $X_j^0$ into $N_j^0$ is a homotopy equivalence, $Y_j'$ is also a relative compact core for $N_j^0$. Since $\widetilde{X}_j$ is precisely invariant under $\Gamma_j$ in $\Gamma$, the restriction of $p_j$ to $Y_j'$ is injective. Hence, $Y_j = p_j(Y_j')$ is a relative compact carrier for $\Gamma_j$.

We have already observed that $Y_j$ and $Y_k$ are disjoint if $j \neq k$ and $1 \leq j, k \leq l$. If $1 \leq j \leq l$ and $l+1 \leq k \leq n$, then since $\widetilde{X}_k$ is disjoint from each translate of $\widetilde{B}_j$, $Y_j$ and $Y_k$ are disjoint. Similarly, if $j \neq k$ and $l+1 \leq j, k \leq n$, we have observed that every translate of $\widetilde{X}_j$ is disjoint from $\widetilde{X}_k$, which implies that $Y_j$ and $Y_k$ are disjoint. Therefore, we have obtained the desired collection of relative compact carriers. □

## 5. Shuffle homotopy equivalences

We introduce in this section the general theory of *shuffle homotopy equivalences*, or *shuffles*, defined to be homotopy equivalences which are homeomorphisms off of specified incompressible submanifolds. We observe that any shuffle has a homotopy inverse which is also a shuffle, and that compositions of shuffles are shuffles.

For $i = 1, 2$, let $M_i$ be a compact, orientable, irreducible 3-manifold and let $\mathcal{V}_i$ be a codimension-zero submanifold of $M_i$. Denote by $\mathrm{Fr}(\mathcal{V}_i)$ the frontier of $\mathcal{V}_i$ in $M_i$. (The frontier is the topological boundary, which equals $\overline{\partial \mathcal{V}_i - \partial M_i}$. We use the term frontier to avoid confusion with the manifold boundary $\partial \mathcal{V}_i$.) We always assume that $\mathrm{Fr}(\mathcal{V}_i)$ is incompressible in $M_i$; in particular, no component of $\mathrm{Fr}(\mathcal{V}_i)$ is a 2-sphere or a boundary-parallel 2-disc. To avoid trivial cases, we assume that $\mathcal{V}_i$ is a nonempty proper subset of $M_i$. Since $\mathrm{Fr}(\mathcal{V}_i)$ consists of a nonempty collection of incompressible surfaces, $\pi_1(M_i)$ is infinite, and each component of $\mathcal{V}_i$ is either a 3-ball or has infinite fundamental group. Thus, since $M_i$ and $\mathcal{V}_i$ are irreducible, both $M_i$ and $\mathcal{V}_i$ are aspherical. A homotopy equivalence $h\colon M_1 \to M_2$ is a *shuffle*, with respect to $\mathcal{V}_1$ and $\mathcal{V}_2$, if $h^{-1}(\mathcal{V}_2) = \mathcal{V}_1$ and $h$ restricts to a homeomorphism from $\overline{M_1 - \mathcal{V}_1}$ to $\overline{M_2 - \mathcal{V}_2}$.

We begin by checking that a shuffle must restrict to a homotopy equivalence from $\mathcal{V}_1$ to $\mathcal{V}_2$.

LEMMA 5.1.   *If $h\colon M_1 \to M_2$ is a shuffle with respect to $\mathcal{V}_1$ and $\mathcal{V}_2$, then the restriction $h\colon (\mathcal{V}_1, \mathrm{Fr}(\mathcal{V}_1)) \to (\mathcal{V}_2, \mathrm{Fr}(\mathcal{V}_2))$ is a homotopy equivalence of pairs.*

*Proof of Lemma* 5.1. We first show that $h$ determines a bijection between the sets of components of $\mathcal{V}_1$ and of $\mathcal{V}_2$. It follows easily from the definition that each component of $\mathcal{V}_2$ contains the image of a component of $\mathcal{V}_1$. Suppose for contradiction that $\mathcal{V}_1$ has more components than $\mathcal{V}_2$. For $i = 1, 2$, construct a graph $\Gamma_i$ with a vertex for each component of $\overline{M_i - \mathcal{V}_i}$ and each component of $\mathcal{V}_i$, and with an edge for each component of the frontier of $\mathcal{V}_i$, which connects the vertices corresponding to the components of $\overline{M_i - \mathcal{V}_i}$ and $\mathcal{V}_i$ that contain it. There is an obvious surjective map $p_i\colon M_i \to \Gamma_i$, which sends a product neighborhood of each component of $\mathrm{Fr}(\mathcal{V}_i)$ onto the corresponding edge of $\Gamma_i$, and the rest of $M_i$ to the vertices of $\Gamma_i$. Note that $p_i$ induces a surjection on fundamental groups. There is also a map $\overline{h}\colon \Gamma_1 \to \Gamma_2$, induced by $h$, so that $\overline{h} \circ p_1$ is homotopic to $p_2 \circ h$. As $h$ is a shuffle, $\overline{M_1 - \mathcal{V}_1}$ and $\overline{M_2 - \mathcal{V}_2}$ have the same number of components, while by assumption $\mathcal{V}_1$ has more components than $\mathcal{V}_2$. Thus $\Gamma_2$ has fewer vertices than $\Gamma_1$ and the same number of edges. Therefore the free group $\pi_1(\Gamma_2)$ has larger rank than $\pi_1(\Gamma_1)$, so that $\overline{h} \circ p_1$ cannot induce a surjection on fundamental groups. But $p_2 \circ h$ does induce a surjection on fundamental groups, a contradiction. Thus, $h$ determines a bijection between the components of $\mathcal{V}_1$ and of $\mathcal{V}_2$.

Let $V_1$ be a component of $\mathcal{V}_1$ and let $V_2$ be the component of $\mathcal{V}_2$ which contains $h(V_1)$. Changing $h$ by a homotopy which respects the hypotheses, we may assume that there is a basepoint $v_1$ in the topological interior of $V_1$ that $h$ maps to a basepoint $v_2$ in the topological interior of $V_2$. Since $h$ is a homotopy equivalence and the frontiers of the $\mathcal{V}_i$ are incompressible, $h$ must induce an injection from $\pi_1(V_1)$ to $\pi_1(V_2)$. To see that $h$ also induces a surjection from $\pi_1(V_1)$ to $\pi_1(V_2)$, let $\beta$ be an element of $\pi_1(V_2)$. Since $h_*:\pi_1(M_1) \to \pi_1(M_2)$ is surjective, there exists a loop $\alpha:(S^1, s_0) \to (M_1, v_1)$ representing an element of $\pi_1(M_1)$ carried to $\beta$ by $h_*$. We may assume that $\alpha$ meets $\mathrm{Fr}(\mathcal{V}_1)$ transversely. There is a map $H$ of $S^1 \times I$ into $M_2$ which maps $S^1 \times \{0\}$ according to $h \circ \alpha$, maps $s_0 \times I$ to $v_2$, and maps $S^1 \times \{1\}$ into $V_2$. We may assume that $H$ is transverse to $\mathrm{Fr}(\mathcal{V}_2)$, and using the incompressibility of $\mathrm{Fr}(\mathcal{V}_2)$, we may assume that the preimage of $\mathrm{Fr}(\mathcal{V}_2)$ consists only of arcs and circles essential in $S^1 \times I$. Since $s_0 \times I$ maps to $v_2$, there are no essential circles, and since $S^1 \times \{1\}$ maps into $V_2$, all the arcs have endpoints in $S^1 \times \{0\}$. Let $\gamma$ be one of these arcs.

There is a disc $E$ in $S^1 \times I$, disjoint from $s_0 \times I$, whose boundary consists of $\gamma$ and a portion $\delta$ of $S^1 \times \{0\}$. Since $h$ is a homeomorphism from $\mathrm{Fr}(\mathcal{V}_1)$ to $\mathrm{Fr}(\mathcal{V}_2)$, there is a unique path $\gamma'$ mapping into $\mathrm{Fr}(\mathcal{V}_1)$ such that $h \circ \gamma'$ is $H|_\gamma$. Replacing $S^1 \times I$ by the closure of $(S^1 \times I) - E$ and $\alpha$ by a new loop with the portion previously represented by $\alpha \circ \delta$ replaced by $\gamma'$, and making a small adjustment, we have a new loop $\alpha$ and a new homotopy $H$ which has fewer arcs in the preimage of $\mathrm{Fr}(\mathcal{V}_2)$. Repeating, we eventually find a loop $\alpha$ in $V_1$ such that $h \circ \alpha$ is homotopic in $V_2$ to $\beta$, showing that $h_*:\pi_1(V_1) \to \pi_1(V_2)$ is surjective. Since $V_1$ and $V_2$ are aspherical, we conclude that the restriction of $h$ to $V_1$ is a homotopy equivalence to $V_2$.

Finally, the restriction $h:(\mathcal{V}_1, \mathrm{Fr}(\mathcal{V}_1)) \to (\mathcal{V}_2, \mathrm{Fr}(\mathcal{V}_2))$ is a homotopy equivalence from $\mathcal{V}_1$ to $\mathcal{V}_2$ and is a homotopy equivalence from $\mathrm{Fr}(\mathcal{V}_1)$ to $\mathrm{Fr}(\mathcal{V}_2)$, since there it is a homeomorphism. It follows from basic theorems of algebraic topology that $h$ is a homotopy equivalence of pairs; see for example Theorem V.3.8 of Whitehead [43]. □

The next lemma is a converse to the previous one. More importantly, it shows that a shuffle has a homotopy inverse which is a shuffle, and agrees with the actual inverse on the complements of the $\mathcal{V}_i$.

LEMMA 5.2. *For $i = 1, 2$, let $M_i$ be a compact, orientable, irreducible 3-manifold, and let $\mathcal{V}_i$ be a codimension-zero submanifold of $M_i$ whose frontier $\mathrm{Fr}(\mathcal{V}_i)$ is nonempty and incompressible. Let $h:M_1 \to M_2$ be a map such that*

(i) $h^{-1}(\mathcal{V}_2) = \mathcal{V}_1$ *and* $h|_{\overline{M_1 - \mathcal{V}_1}}: \overline{M_1 - \mathcal{V}_1} \to \overline{M_2 - \mathcal{V}_2}$ *is a homeomorphism, and*

(ii) $h|_{\mathcal{V}_1}\colon \mathcal{V}_1 \to \mathcal{V}_2$ *is a homotopy equivalence.*

*Then, $h$ is a homotopy equivalence, and there exists a homotopy inverse $\overline{h}\colon M_2 \to M_1$ for $h$ so that $\overline{h}^{-1}(\mathcal{V}_1) = \mathcal{V}_2$, so that $\overline{h}|_{\overline{M_2-\mathcal{V}_2}}$ is the inverse of $h|_{\overline{M_1-\mathcal{V}_1}}$, and so that $\overline{h}|_{\mathcal{V}_2}\colon \mathcal{V}_2 \to \mathcal{V}_1$ is a homotopy equivalence. Moreover, $\overline{h}\circ h$ is homotopic to the identity relative to $\overline{M_1-\mathcal{V}_1}$, and $h\circ\overline{h}$ is homotopic to the identity relative to $\overline{M_2-\mathcal{V}_2}$.*

*Proof of Lemma* 5.2. The restriction of $h$ to $\mathcal{V}_1$ defines a map of pairs $h_0\colon(\mathcal{V}_1,\mathrm{Fr}(\mathcal{V}_1)) \to (\mathcal{V}_2,\mathrm{Fr}(\mathcal{V}_2))$ which is a homeomorphism from $\mathrm{Fr}(\mathcal{V}_1)$ to $\mathrm{Fr}(\mathcal{V}_2)$. As in the proof of Lemma 5.1, it follows that $h_0$ is a homotopy equivalence of pairs. Let $\overline{h_0}\colon(\mathcal{V}_2,\mathrm{Fr}(\mathcal{V}_2)) \to (\mathcal{V}_1,\mathrm{Fr}(\mathcal{V}_1))$ be a homotopy inverse for $h_0$. On $\mathrm{Fr}(\mathcal{V}_1)$, $\overline{h_0}\circ h_0$ is homotopic to the identity, so that $\overline{h_0}|_{\mathrm{Fr}(\mathcal{V}_2)}$ is homotopic to $\left(h_0|_{\mathrm{Fr}(\mathcal{V}_1)}\right)^{-1}$. So we may assume that $\overline{h_0}|_{\mathrm{Fr}(\mathcal{V}_2)} = \left(h_0|_{\mathrm{Fr}(\mathcal{V}_1)}\right)^{-1}$.

Consider next a homotopy of pairs $K$ from $\overline{h_0}\circ h_0$ to the identity map $\mathrm{id}_{(\mathcal{V}_1,\mathrm{Fr}(\mathcal{V}_1))}$. On a 2-manifold, any homotopy from the identity to the identity can be deformed to an isotopy from the identity to the identity; so, using the homotopy extension property, we may assume that $K$ restricts on $\mathrm{Fr}(\mathcal{V}_1)\times I$ to an isotopy $J$ from $\mathrm{id}_{\mathrm{Fr}(\mathcal{V}_1)}$ to $\mathrm{id}_{\mathrm{Fr}(\mathcal{V}_1)}$. Consequently, we may change $\overline{h_0}$ by a homeomorphism supported in a collar neighborhood of $\mathrm{Fr}(\mathcal{V}_1)$ in $\mathcal{V}_1$ to ensure that $\overline{h_0}\circ h_0$ is homotopic to $\mathrm{id}_{\mathcal{V}_1}$ relative to $\mathrm{Fr}(\mathcal{V}_1)$. (Explicitly, let $\mathrm{Fr}(\mathcal{V}_1)\times[0,1]$ be a collar neighborhood of $\mathrm{Fr}(\mathcal{V}_1)$ in $\mathcal{V}_1$ with $\mathrm{Fr}(\mathcal{V}_1) = \mathrm{Fr}(\mathcal{V}_1)\times\{0\}$, and define the homeomorphism $j$ of $\mathcal{V}_1$ by $j(x,s) = (J_s^{-1}(x),s)$ in the collar and $j(v) = v$ for $v$ not in the collar. An isotopy $G$ from $j$ to $\mathrm{id}_{\mathcal{V}_1}$ is defined by letting $G_t(x,s) = (J_{t+s}^{-1}(x),s)$ for $0 \leq s \leq 1-t$ and $G_t(v) = v$ for all other $v \in \mathcal{V}_1$. Since $(G_t \circ K_t)(x,0) = G_t(J_t(x),0) = (x,0)$ for all $t$, the map $L\colon(\mathcal{V}_1\times I, \mathrm{Fr}(\mathcal{V}_1)\times I) \to (\mathcal{V}_1,\mathrm{Fr}(\mathcal{V}_1))$ given by letting $L_t = G_t\circ K_t$ for all $t$ is a homotopy from $(j\circ\overline{h_0})\circ h_0$ to $\mathrm{id}_{\mathcal{V}_1}$ relative to $\mathrm{Fr}(\mathcal{V}_1)$.)

Similarly, there exists a left homotopy inverse $h_1$ of $\overline{h_0}$ so that $h_1\circ\overline{h_0}$ is homotopic to $\mathrm{id}_{\mathcal{V}_2}$ relative to $\mathrm{Fr}(\mathcal{V}_2)$. Then, we have $h_1 \simeq h_1\circ\overline{h_0}\circ h_0 \simeq h_0$, with all homotopies relative to $\mathrm{Fr}(\mathcal{V}_1)$, and so $h_0\circ\overline{h_0} \simeq \mathrm{id}_{\mathcal{V}_2}$ relative to $\mathrm{Fr}(\mathcal{V}_2)$. Now we can define $\overline{h}\colon M_2\to M_1$ by using $\overline{h_0}$ on $\mathcal{V}_2$ and the inverse of $h|_{\overline{M_1-\mathcal{V}_1}}$ on $\overline{M_2-\mathcal{V}_2}$. One easily checks that $\overline{h}$ has the desired properties. □

## 6. The characteristic submanifold

In this section we prove that when one performs a shuffle with respect to the characteristic submanifold of a hyperbolizable 3-manifold, the image of the characteristic submanifold of the domain is the characteristic submanifold of the range. We begin by recalling the definition of the characteristic submanifold.

If $X \subseteq Y$, then by $\pi_1(X) \to \pi_1(Y)$ we mean the homomorphism induced by inclusion. When this homomorphism is injective, we regard $\pi_1(X)$ as a subgroup of $\pi_1(Y)$.

A surface $X$ is *properly embedded* in a 3-manifold $M$ if $\partial X = X \cap \partial M$. A torus or properly embedded annulus $X$ is *essential* if $\pi_1(X) \to \pi_1(M)$ is injective and $X$ is not properly homotopic into $\partial M$. If $M$ is compact, orientable, irreducible and contains no essential tori, then it is said to be *atoroidal*. A properly embedded annulus $X$ in $M$ is *primitive* if it is incompressible and the generator of $\pi_1(X)$ is not a proper power in $\pi_1(M)$.

An embedded $I$-bundle $R$ in $M$ is *admissibly embedded* if $R \cap \partial M$ is the associated $\partial I$-bundle of $R$. An embedded Seifert-fibered space $R$ in $M$ is *admissibly embedded* if $R \cap \partial M$ is a union of fibers in $\partial R$. An admissibly embedded $I$-bundle or Seifert-fibered space $R$ in $M$ is *essential* if every component of the frontier of $R$ in $M$ is an essential torus or annulus in $M$. In particular, this implies that $\pi_1(R) \to \pi_1(M)$ is injective. A homotopy $F\colon R \times I \to M$ is *admissible* if $F_t(F_0^{-1}(\partial M)) \subseteq \partial M$ for all $t$.

A compact codimension-zero submanifold $\Sigma$ of $M$ has the *engulfing property* if every essential embedding $f\colon R \to M$ of an $I$-bundle or a Seifert-fibered space into $M$ is admissibly homotopic to a map with image in $\Sigma$. We define $\Sigma$ to be a *characteristic submanifold* of $M$ if $\Sigma$ consists of a collection of essential $I$-bundles and Seifert-fiber spaces having the engulfing property, and $\Sigma$ is minimal in the sense that no proper subcollection of the components of $\Sigma$ has the engulfing property. Jaco and Shalen [18] and Johannson [19] show that every compact, orientable, irreducible 3-manifold $M$ with nonempty incompressible boundary contains a characteristic submanifold, and that any two characteristic submanifolds are admissibly isotopic in $M$. Hence, we often speak of *the* characteristic submanifold $\Sigma(M)$ of $M$.

By the Squeezing Theorem in Section V.1 of Jaco-Shalen [18] or by Proposition 10.8 of Johannson [19], the characteristic submanifold satisfies a stronger engulfing property: every essentially embedded fibered submanifold is admissibly isotopic into $\Sigma$.

The characteristic submanifold we use in this paper is that of Jaco and Shalen. It differs slightly from the one given in Johannson's formulation (using the boundary pattern consisting of the set of boundary components of $M$). The discrepancy arises from the fact that in Johannson's theory, any incompressible torus is considered to be essential, even if it is homotopic into $\partial M$. In particular, Johannson's characteristic submanifold must contain every torus boundary component of $M$, so it must be allowed to have components that are regular neighborhoods of torus boundary components. If such components are deleted, one obtains the characteristic submanifold of Jaco and Shalen.

When $M$ is a compact, hyperbolizable 3-manifold with nonempty incompressible boundary, every Seifert-fibered component $V$ of $\Sigma(M)$ is homeomor-phic to either a solid torus or a thickened torus; see Morgan [33] or Canary and McCullough [12]. If $V$ is a solid torus, its frontier in $M$ consists of a nonempty collection of essential annuli in $M$. If $V$ is a thickened torus, one of its boundary tori is a boundary component of $M$, and its other boundary torus meets $\partial M$ in a nonempty collection of annuli which are pairwise nonhomotopic in $\partial M$. Note that every component of the frontier of $\Sigma(M)$ is an essential annulus. If $\Sigma(M)$ contains a component $V$ homeomorphic to a thickened torus such that at least two components of $V \cap \partial M$ are annuli, then $M$ is said to have *double trouble*.

We are now prepared to prove the main result of the section.

PROPOSITION 6.1.   *Let $M_1$ and $M_2$ be compact, hyperbolizable 3-manifolds with nonempty incompressible boundary, and let $\mathcal{V}_1$ and $\mathcal{V}_2$ be essential fibered submanifolds of $M_1$ and $M_2$. Let $f\colon M_1 \to M_2$ be a shuffle with respect to $\mathcal{V}_1$ and $\mathcal{V}_2$. Then $\mathcal{V}_1$ is a characteristic submanifold for $M_1$ if and only if $\mathcal{V}_2$ is a characteristic submanifold for $M_2$.*

*Proof of Proposition* 6.1. Assume that $\mathcal{V}_2$ is a characteristic submanifold for $M_2$. Since $M_2$ is hyperbolizable, every component of $\mathrm{Fr}(\mathcal{V}_2)$ is an annulus. Since $f|_{\overline{M_1-\mathcal{V}_1}}$ is a homeomorphism, every component of $\mathrm{Fr}(\mathcal{V}_1)$ is an annulus as well.

Let $\Sigma$ be a characteristic submanifold for $M_1$. Then $\mathcal{V}_1$ is isotopic into $\Sigma$, so that we may choose $\Sigma$ so that $\mathcal{V}_1$ lies in the topological interior of $\Sigma$. Since $\mathrm{Fr}(\mathcal{V}_1)$ consists of annuli, Corollary 5.7 of Johannson [19] shows that $\Sigma$ admits an admissible fibering such that $\mathrm{Fr}(\mathcal{V}_1)$ is vertical (i.e. is a union of fibers). This shows that $\overline{\Sigma - \mathcal{V}_1}$ is an essential fibered submanifold of $M_1$, and moreover that each component of $\mathcal{V}_1$ is an admissibly fibered $I$-bundle, solid torus, or thickened torus (since these are the possible components for $\Sigma$ and $\mathcal{V}_1$ is an essential fibered submanifold of $\Sigma$). We will prove that each component $R_1$ of $\overline{\Sigma - \mathcal{V}_1}$ is a product region with one end a component of $\mathrm{Fr}(\mathcal{V}_1)$ and the other end a component of $\mathrm{Fr}(\Sigma)$ disjoint from $\mathcal{V}_1$. This implies that $\mathcal{V}_1$ is admissibly isotopic to $\Sigma$, and hence is characteristic.

Let $R_2$ be the image of $R_1$ under the homeomorphism $f|_{\overline{M_1-\mathcal{V}_1}}$. Since $R_1$ is admissibly fibered, so is $R_2$. We claim that $R_2$ is essential. Since $\mathcal{V}_2$ is essential, $R_2 \cap \mathcal{V}_2$ is essential. Suppose that $F_2$ is a component of $\mathrm{Fr}(R_2)$ that does not lie in $\mathcal{V}_2$. If $F_2$ is not essential, there is a proper homotopy carrying $F_2$ into $\partial M_2$. Since $\mathrm{Fr}(\mathcal{V}_2)$ is essential, this homotopy may be deformed off of $\mathrm{Fr}(\mathcal{V}_2)$ to give a proper homotopy of $F_2$ into $\partial M_2$ which moves $F_2$ only through $\overline{M_2 - \mathcal{V}_2}$. Under $\left(f|_{\overline{M_1-\mathcal{V}_1}}\right)^{-1}$, this corresponds to a proper homotopy of a component of $\mathrm{Fr}(\Sigma)$ into $\partial M_1$, a contradiction since $\Sigma$ is essential. Therefore $F_2$ is essential, and the claim is proved.

By the engulfing property of $\mathcal{V}_2$, $R_2$ is admissibly homotopic into $\mathcal{V}_2$. For later reference, we isolate the next portion of the argument as a lemma.

LEMMA 6.2. *Let $M$ be a compact, orientable, irreducible 3-manifold. Let $R$ be a compact connected 3-manifold embedded in $M$, with incompressible frontier, which is homotopic into $M - R$. Then $R$ is a product of the form $F \times I$ where $F \times \{0\}$ is a component of $\mathrm{Fr}(R)$. If the homotopy is admissible and every annular component of $\mathrm{Fr}(R)$ is essential, then the product structure can be chosen so that $\mathrm{Fr}(R)$ equals $F \times \partial I$.*

*Proof of Lemma* 6.2. Since $\mathrm{Fr}(R)$ is incompressible, $R$ is also irreducible. Consider the covering space $\widehat{M}$ of $M$ corresponding to the subgroup $\pi_1(R)$ of $\pi_1(M)$. It contains a lift $\widehat{R}$ of $R$, and $\pi_1(\widehat{R}) \to \pi_1(\widehat{M})$ is an isomorphism. This implies that each component $W$ of the closure of $\widehat{M} - \widehat{R}$ meets one component $G_W$ of $\mathrm{Fr}(\widehat{R})$. Since $\mathrm{Fr}(R)$ is incompressible, so is $\mathrm{Fr}(\widehat{R})$; hence $\pi_1(G_W) \to \pi_1(W)$ is an isomorphism. Let $H$ be a homotopy carrying $R$ into $M - R$. It lifts to a homotopy of $\widehat{R}$ that moves $\widehat{R}$ into one of the components $X$ of the closure of $\widehat{M} - \widehat{R}$. This implies that the subgroup $\pi_1(X)$ equals $\pi_1(\widehat{R})$, and hence that $\pi_1(G_X) \to \pi_1(\widehat{R})$ is an isomorphism. By Theorem 10.5 of [15], $\widehat{R}$ is a product $G_X \times I$ with $G_X = G_X \times \{0\}$.

Assume now that $H$, and hence its lift to $\widehat{M}$, is admissible. Suppose that $A$ is a component of $\widehat{R} \cap \partial \widehat{M}$. Since $\mathrm{Fr}(\widehat{R})$ is incompressible, the components of $\partial A$ are not contractible in $\partial \widehat{M}$. The lifted homotopy carries $A$ into $\partial \widehat{M} - A$. By the 2-dimensional analogue of the previous paragraph, $A$ is an annulus, whose boundary curves are parallel in $\partial \widehat{M}$ to a component $C$ of $\partial G_X$.

Since all components of $\widehat{R} \cap \partial \widehat{M}$ are annuli, we may choose the product structure $G_X \times I$ of $\widehat{R}$ so that the components of $\widehat{R} \cap \partial \widehat{M}$ that meet $G_X$ form $\partial G_X \times I$. To show that $\widehat{R}$, and hence $R$, have the desired product structures, we must prove that all components of $\widehat{R} \cap \partial \widehat{M}$ meet $\partial G_X$.

To motivate the remainder of the argument, we first give an example showing how components of $\widehat{R} \cap \partial \widehat{M}$ can fail to meet $\partial G_X$ when $\mathrm{Fr}(\widehat{R})$ contains inessential annuli. Start with $G \times [0,1] \subset G \times \mathbf{R}$, where $G$ is a compact surface with nonempty boundary. Let $C$ be a boundary component of $G$, and consider a small regular neighborhood of $C \times \{1\}$ in $G \times [1, \infty)$. The frontier of this neighborhood in $G \times [1, \infty)$ is a properly embedded annulus; let $N$ be a small regular neighborhood of this annulus in $G \times [1, \infty)$. Then $G \times [0,1] \cup N$ is admissibly homotopic into $G \times (-\infty, 0)$, and its frontier consists of two copies of $G$ plus an inessential annulus.

Assume now that every annular component of $\mathrm{Fr}(R)$, and hence of $\mathrm{Fr}(\widehat{R})$, is essential, and suppose for contradiction that some annulus $A$ of $\widehat{R} \cap \partial \widehat{M}$ does not meet $G_X$. Since the boundary components of $A$ are parallel to some component $C$ of $\partial G_X$, there is an annular component $B$ of $\partial \widehat{M} - \widehat{R}$. Let $Y$ be the component of the closure of $\widehat{M} - \widehat{R}$ that contains $B$. Since $G_Y$ is the entire frontier of $Y$, $G_Y \cup B$ must be a component of $\partial Y$. Let $K$ be a compact core of $Y$ that contains $G_Y \cup B$. Since $\pi_1(G_Y) \to \pi_1(K)$ is an isomorphism, Theorem 10.5 of [15] shows that $K$ is a product $I$-bundle with one component of the associated $\partial I$-bundle equal to $G_Y$. Since the other component must lie in $B$, $G_Y$ must be an annulus. Therefore $K$ is a solid torus in which $G_Y$ is parallel to $B$. But then, $G_Y$ is inessential. This contradiction proves that $\widehat{R} \cap \partial \widehat{M} = \partial G_X \times I$. □

Recall that $R_2$ is the image of $R_1$ under the homeomorphism $f|_{\overline{M_1 - \mathcal{V}_1}}$. Lemma 6.2 implies that $R_2$ and hence $R_1$ have product structures whose ends equal the components of their frontiers, all of which are annuli. We now show that exactly one component of $\mathrm{Fr}(R_2)$ lies in $\mathcal{V}_2$.

Suppose first that $R_2$ is disjoint from $\mathcal{V}_2$, and hence that $R_1$ is a component of $\Sigma$ which does not meet $\mathcal{V}_1$. Let $F$ be a component of the frontier of $R_2$. Using Proposition 19.1 of Johannson [19], one may show that, since $F$ is properly homotopic into $\mathcal{V}_2$, it must be parallel in $\overline{M_2 - \mathcal{V}_2}$ to a component of $\mathrm{Fr}(\mathcal{V}_2)$. Since $f|_{\overline{M_1 - \mathcal{V}_1}}$ is a homeomorphism, the components of $\mathrm{Fr}(R_1)$ are parallel to components of $\mathrm{Fr}(\mathcal{V}_1)$. It follows that $R_1$ is admissibly homotopic into $\mathcal{V}_1$. Therefore $\Sigma - R_1$ still has the engulfing property, contradicting the fact that the characteristic submanifold is a minimal submanifold having the engulfing property.

Suppose for contradiction that both components of $\mathrm{Fr}(R_2)$ lie in $\mathcal{V}_2$. Then both components of $\mathrm{Fr}(R_1)$ lie in $\mathcal{V}_1$. Since $\Sigma$ can be fibered so that $\mathcal{V}_1$ is a union of fibers, the components (not necessarily distinct) of $\mathcal{V}_1$ that meet $R_1$ can both be $I$-fibered or can both be Seifert-fibered. Therefore the same is true for the components of $\mathcal{V}_2$ that meet $R_2$. (A component of $\mathcal{V}_1$ can be $I$-fibered exactly when its fundamental group is free of rank at least 2, or when it is free of rank 1 and the component is admissibly fibered as an $I$-bundle over the annulus or Möbius band, and these properties are preserved by the homotopy equivalence of pairs $(\mathcal{V}_1, \mathrm{Fr}(\mathcal{V}_1)) \to (\mathcal{V}_2, \mathrm{Fr}(\mathcal{V}_2))$. A component of $\mathcal{V}_1$ can be Seifert-fibered if and only if its fundamental group contains an infinite cyclic normal subgroup.) Since up to isotopy, the annulus admits a unique $I$-fibering and a unique $S^1$-fibering, and $R_2$ is the product of an annulus and an interval, the fiberings on these components of $\mathcal{V}_2$ extend over $R_2$. Therefore $\mathcal{V}_2 \cup R_2$ can be admissibly fibered. Now $\mathcal{V}_2 \cup R_2$ has the engulfing property, since $\mathcal{V}_2$ does. But $\mathcal{V}_2 \cup R_2$ is not homeomorphic to $\mathcal{V}_2$, so is not characteristic. Hence, the union of some proper subcollection of the components of $\mathcal{V}_2$ has the engulfing property. A minimal such union would give a characteristic submanifold with fewer components than $\mathcal{V}_2$, contradicting the uniqueness of $\mathcal{V}_2$.

We have shown that exactly one component of $\mathrm{Fr}(R_2)$ lies in $\mathcal{V}_2$. It follows that $\mathrm{Fr}(R_1)$ has exactly one component in $\mathcal{V}_1$. The other must be a component of Fr($\Sigma$). Since this is true for all components of $\overline{\Sigma - \mathcal{V}_1}$, we conclude that $\Sigma$ is isotopic to $\mathcal{V}_1$, hence that $\mathcal{V}_1$ is characteristic.

Conversely, suppose that $\mathcal{V}_1$ is characteristic. Let $g$ be a homotopy inverse of $f$ obtained using Lemma 5.2. The previous case applies with $g$ in the role of $f$, showing that $\mathcal{V}_2$ is characteristic. □

## 7. Primitive shuffle equivalence

We now begin to specialize the discussion of the previous several sections to the setting of shuffles with respect to collections of primitive solid tori. A solid torus $V$ embedded in a compact, orientable, irreducible 3-manifold is *primitive* if its frontier consists of primitive annuli. In particular, the inclusion of each component of the frontier of $V$ into $V$ is a homotopy equivalence. In Lemma 7.1 we show that two shuffles which agree off of collections of primitive solid tori must differ up to homotopy by Dehn twists about frontier annuli. We then formally define primitive shuffle equivalence and apply Lemma 7.1 to show that primitive shuffle equivalence determines a finite-to-one equivalence relation on $\mathcal{A}(M)$.

Let $M$ be a compact, orientable, irreducible 3-manifold, and let $A$ be an embedded annulus in $M$ with $A \cap \partial M = \partial A$. Let $A \times [0, 1]$ be a collar neighborhood on one side of $A$, meeting $\partial M$ in $\partial A \times [0, 1]$. Regarding $A$ as $S^1 \times [0, 1]$, let $h_n \colon A \times [0, 1] \to A \times [0, 1]$ be a homeomorphism defined by $h_n((\exp(2\pi i u), v), w) = ((\exp(2\pi i (u + nw)), v), w)$. Thus $h_n$ is the identity on $A \times \{0, 1\}$, and the annuli $A \times \{w\}$ rotate $n$ full turns as one moves across the $[0, 1]$-factor of $A \times [0, 1]$. A *Dehn twist about $A$* is a homeomorphism of $M$ which equals $h_n$ on $A \times [0, 1]$ and is the identity map on the complement of $A \times [0, 1]$. The following lemma relates Dehn twists about annuli to shuffles.

LEMMA 7.1. *Let $s_0, s_1 \colon M_1 \to M_2$ be shuffles with respect to $\mathcal{V}_i$, where the $\mathcal{V}_i$ are unions of finitely many disjoint primitive solid tori in $M_i$. If $s_0|_{\overline{M_1 - \mathcal{V}_1}} = s_1|_{\overline{M_1 - \mathcal{V}_1}}$, then*

(i) *there is a homeomorphism $r \colon M_1 \to M_1$, which is a composition of Dehn twists about frontier annuli of $\mathcal{V}_1$, such that $s_0 \circ r$ is homotopic to $s_1$ relative to $\overline{M_1 - \mathcal{V}_1}$, and*

(ii) *there is a homeomorphism $r' \colon M_2 \to M_2$, which is a composition of Dehn twists about frontier annuli of $\mathcal{V}_2$, such that $r' \circ s_0$ is homotopic to $s_1$ relative to $\overline{M_1 - \mathcal{V}_1}$.*

*Proof of Lemma* 7.1. Let $V_1$ be a component of $\mathcal{V}_1$, and let $V_2$ be the component of $\mathcal{V}_2$ that contains $s_j(V_1)$. Let $r_j = s_j|_{V_1}$. Fix a frontier annulus $A$ of $V_1$, and let $C$ be its core curve. Since the $r_j$ agree on $A$, they induce the same isomorphism from $\pi_1(V_1)$ to $\pi_1(V_2)$, and therefore are homotopic. Let $H\colon V_1 \times I \to V_2$, be a homotopy from $r_0$ to $r_1$ and let $h = H|_{C \times I}$. Since $r_0$ and $r_1$ agree on $C$, $h$ induces a map $\bar{h}$ from the torus $W = C \times I/\langle (x,0) \sim (x,1)\rangle$ into $V_2$. Let $C'$ be a simple loop in $W$ that intersects $C$ in one point. Since $A$ is primitive, $C$ represents a generator of $\pi_1(V_1)$. Therefore we can write $\bar{h}_*(C') = n \cdot C$ for some $n$. Then $\bar{h}_*(C' - n \cdot C) = 0$ in $\pi_1(V_2)$. Note also that $C' - n \cdot C$ is representable by a simple loop in $W$ that meets $C$ in one point. This implies that $h\colon C \times I \to V_2$ is homotopic, relative to $C \times \{0,1\}$, to a map into $C$. (To see this, note that, since $\bar{h}_*(C' - n \cdot C) = 0$, $\bar{h}$ extends to the union of $W$ with a 2-disc attached along $C' - n \cdot C$. Attaching a 3-ball to this union, one obtains a solid torus $V$, and $\bar{h}$ extends over $V$ since $\pi_2(V_2) = 0$. There is a homotopy of $W$, relative to $C$, that moves it through $V$ onto $C$. This homotopy then gives the desired homotopy of $h$.)

By the homotopy extension property, we may assume that $h_t(C) = r_0(C)$ for all $t$. By a further adjustment, we may assume that $H_t(A) = r_0(A)$ for all $t$. Any homotopy from $\mathrm{id}_A$ to $\mathrm{id}_A$ can be deformed to an isotopy $J$ that rotates $A$ some number of times in the $S^1$-direction. So we may assume that $H_t$ has the form $r_0 \circ J_t$ on $A$ for all $t$. If $r_0$ is precomposed by the correct Dehn twist about $A$, the composition will be homotopic to $r_1$ relative to $A$. (A detailed construction of this type is given in the proof of Lemma 5.2.) Applying this reasoning to each of the frontier annuli of $\mathcal{V}_1$, we obtain a product of Dehn twists $r$ such that $r_0 \circ r$ is homotopic to $r_1$ relative to the frontier of $\mathcal{V}_1$. Extending this homotopy to $M_1$ using the homeomorphism $s_0|_{\overline{M_1 - \mathcal{V}_1}}$ at each level gives the homotopy from $s_0 \circ r$ to $s_1$. This proves (i).

For (ii), let $s_0'$ and $s_1'$ be the homotopy inverses for $s_0$ and $s_1$ obtained using Lemma 5.2. By (i), there exists $r'\colon M_2 \to M_2$ such that $s_1' \circ r' \simeq s_0'$ relative to $\overline{M_2 - \mathcal{V}_2}$. Therefore $s_1 \circ s_1' \circ r' \circ s_0 \simeq s_1 \circ s_0' \circ s_0$ relative to $\overline{M_1 - \mathcal{V}_1}$, and so $r' \circ s_0 \simeq s_1$ relative to $\overline{M_1 - \mathcal{V}_1}$. $\square$

Let $M_1$ and $M_2$ be compact, oriented, irreducible 3-manifolds with nonempty incompressible boundary. A shuffle $s\colon M_1 \to M_2$ with respect to $\mathcal{V}_1$ and $\mathcal{V}_2$ is called *primitive* if

1. each $\mathcal{V}_i$ is a collection of primitive solid torus components of $\Sigma(M_i)$, and

2. $s$ restricts to an *orientation-preserving* homeomorphism from $\overline{M_1 - \mathcal{V}_1}$ to $\overline{M_2 - \mathcal{V}_2}$.

If $M$ is a compact, hyperbolizable 3-manifold with nonempty incompressible boundary, two elements $[(M_1, h_1)]$ and $[(M_2, h_2)]$ of $\mathcal{A}(M)$ are *primitive shuffle equivalent* if there exists a primitive shuffle $s\colon M_1 \to M_2$ such that $[(M_2, h_2)] = [(M_2, s \circ h_1)]$. Note that by Lemma 5.2, a primitive shuffle with respect to $\mathcal{V}_1$ and $\mathcal{V}_2$ has a homotopy inverse that is a primitive shuffle with respect to $\mathcal{V}_2$ and $\mathcal{V}_1$. Also, a composition of primitive shuffles is homotopic to a primitive shuffle. For suppose that $s_1\colon M_1 \to M_2$ and $s_2\colon M_2 \to M_3$ are primitive shuffles. If the characteristic submanifold $\Sigma(M_2)$ used in the definition of $s_1$ equals the one used in the definition of $s_2$, then the composition $s_2 \circ s_1$ is a primitive shuffle. If not, then since $\Sigma(M_2)$ is well-defined up to ambient isotopy, $s_2 \circ s_1$ is still homotopic to a primitive shuffle. Therefore, primitive shuffle equivalence determines an equivalence relation on $\mathcal{A}(M)$.

Define $\widehat{\mathcal{A}}(M)$ to be the collection of equivalence classes, and let $q\colon \mathcal{A}(M) \to \widehat{\mathcal{A}}(M)$ be the quotient map. The next result shows that $q$ is finite-to-one.

PROPOSITION 7.2. *Let $M$ be a compact, hyperbolizable 3-manifold with nonempty incompressible boundary. Then the quotient map $q\colon \mathcal{A}(M) \to \widehat{\mathcal{A}}(M)$ is finite-to-one.*

*Proof of Proposition* 7.2. Let $[(M_1, h_1)]$ be a fixed element of $\mathcal{A}(M)$ and suppose that $\mathcal{A}(M)$ contains infinitely many distinct marked homeomorphism types

$$\mathcal{C} = \{[(M_2, h_2)], \ldots, [(M_j, h_j)], \ldots\}$$

which are primitive shuffle equivalent to $[(M_1, h_1)]$. Fix primitive shuffles $s_j\colon M_1 \to M_j$ such that $s_j \circ h_1 \simeq h_j$. For each $j$, let $\overline{s_j}\colon M_j \to M_1$ be a primitive shuffle which is a homotopy inverse of $s_j$.

Let $W_j$ be the union of the components of $\Sigma(M_j)$ that are primitive solid tori, and let $A_j$ be its collection of frontier annuli. By Lemma 5.1, $s_j\colon (W_1, A_1) \to (W_j, A_j)$ is a homotopy equivalence of pairs. Since the annuli are primitive, this implies that each pair $(W_j, A_j)$ is homeomorphic to $(W_1, A_1)$.

Fix homeomorphisms $F_j\colon (W_1, A_1) \to (W_j, A_j)$. Let $f_j$ be the restriction of $F_j$ to a homeomorphism from $A_1$ to $A_j$, and let $g_j\colon A_1 \to A_j$ be the restriction of $s_j$ to $A_1$. Notice that both $f_j$ and $g_j$ are homeomorphisms. Since the mapping class group of $A_1$ is finite, there exist two distinct elements $[(M_k, h_k)]$ and $[(M_l, h_l)]$ of $\mathcal{C}$, such that $f_k^{-1} \circ g_k$ is isotopic to $f_l^{-1} \circ g_l$. Thus, $g_l \circ g_k^{-1}$ is isotopic to $f_l \circ f_k^{-1}$. Changing $F_l$ by an isotopy, we may assume that $g_l \circ g_k^{-1} = f_l \circ f_k^{-1}$.

Define $s' = s_l \circ \overline{s_k}$. Then $s' \circ h_k = s_l \circ \overline{s_k} \circ h_k \simeq s_l \circ \overline{s_k} \circ s_k \circ h_1 \simeq h_l$. Notice that $F_l \circ F_k^{-1}$ is a homeomorphism from $(W_k, A_k)$ to $(W_l, A_l)$ which agrees with $s'$ on $A_k$. Define $s''\colon M_k \to M_l$ using $s'$ on $\overline{M_k - W_k}$ and $F_l \circ F_k^{-1}$ on $W_k$. Lemma 7.1 then implies that there exists an orientation-preserving homeomorphism $r\colon M_k \to M_k$ such that $s'' \circ r$ is homotopic to $s'$. Since $s'' \circ r$ is an orientation-preserving homeomorphism and $s'' \circ r \circ h_k \simeq h_l$, we see that $[(M_k, h_k)] = [(M_l, h_l)]$, which is a contradiction. □

## 8. The continuity of $\widehat{\Psi}$

This section is devoted to the proof of the following more general version of Theorem A. In the present statement $\Psi : \mathcal{D}(\pi_1(M)) \to \mathcal{A}(M)$ denotes the lift of the map $\Psi : AH(\pi_1(M)) \to \mathcal{A}(M)$ defined in Section 2.

THEOREM A. *Let $M$ be a compact, hyperbolizable 3-manifold with incompressible boundary. Suppose that $\{\rho_i\}$ converges to $\rho$ in $\mathcal{D}(\pi_1(M))$. Then, for all sufficiently large $i$, $\Psi(\rho_i)$ is primitive shuffle equivalent to $\Psi(\rho)$. Moreover, if $\{\rho_i(\pi_1(M))\}$ is geometrically convergent (not necessarily to $\rho(\pi_1(M))$), then $\Psi(\rho_i)$ is eventually constant.*

As the argument is rather intricate, we will begin with an outline.

*Outline of the proof of Theorem* A. We first reduce to the case that $\{\rho_i(\pi_1(M))\}$ converges geometrically to a Kleinian group $\widehat{\Gamma}$. Lemma 3.5 produces a sequence $\{f_i \colon Z_i \to N_i\}$ of biLipschitz diffeomorphisms of submanifolds $\{Z_i\}$ of the geometric limit $\widehat{N} = \mathbf{H}^3/\widehat{\Gamma}$ to submanifolds of the approximates $\{N_i = \mathbf{H}^3/\rho_i(\pi_1(M))\}$. If $\widehat{\Gamma} = \rho(\pi_1(M))$ (i. e. the algebraic and geometric limits agree) and $M_0$ is a compact core for the algebraic limit $N = \mathbf{H}^3/\rho(\pi_1(M))$, then $f_i(M_0)$ is a compact core of $N_i$ for all large enough $i$. It is then easy to check that the approximates and the algebraic limit have the same marked homeomorphism type for all large enough $i$. One may apply the same argument whenever there is a compact core for the algebraic limit which embeds (via the obvious covering map $\pi \colon N \to \widehat{N}$) in the geometric limit.

In the general case, we build a compact core $M_1$ for the algebraic limit with the property that after one removes a finite collection $\mathcal{U}$ of solid tori (all lying within rank one cusps), the results of Section 4 may be used to show that the remainder $M_1 - \mathcal{U}$ embeds in the geometric limit $\widehat{N}$ (via the covering map $\pi$). One then constructs a submanifold $M_2$ of the geometric limit $\widehat{N}$, which is homotopy equivalent to $M_1$, by appending a finite collection of solid tori (lying in cusps of the geometric limit) to the image $\pi(M_1 - \mathcal{U})$ of $M_1 - \mathcal{U}$. Moreover, we construct a primitive shuffle $\phi \colon M_1 \to M_2$ which agrees with the covering map wherever possible. We then show that $f_i(M_2)$ is a compact core for $N_i$ for all large enough $i$. One then checks that if $\Psi(\rho) = [(M_1, \psi)]$, then $\Psi(\rho_i) = [(M_2, \phi \circ \psi)]$ for all large enough $i$, which establishes the theorem.

*Proof of Theorem* A. We may assume that $\pi_1(M)$ is nonabelian, since if $\pi_1(M)$ is abelian, then $M$ is a thickened torus and $\mathcal{A}(M)$ contains a single element.

We will show that if $\{\rho_i(\pi_1(M))\}$ converges geometrically to a Kleinian group $\widehat{\Gamma}$, then $\Psi(\rho_i)$ is eventually constant and $\Psi(\rho_i)$ is primitive shuffle equivalent to $\Psi(\rho)$ for all large enough $i$. As every subsequence of $\{\rho_i\}$ has a further subsequence which converges geometrically, this shows that only finitely many of the $\Psi(\rho_i)$ can fail to be primitive shuffle equivalent to $\Psi(\rho)$, proving Theorem A.

Let $N = \mathbf{H}^3/\rho(\pi_1(M))$ and $\widehat{N} = \mathbf{H}^3/\widehat{\Gamma}$, and let $\pi\colon N \to \widehat{N}$ be the covering map. Since $M$ has incompressible boundary, $\pi_1(M)$ is freely indecomposable, so Bonahon's theorem [8] implies that $N$ is topologically tame. Let $\widehat{\mathcal{H}}$ be a precisely invariant system of horoballs for $\widehat{\Gamma}$ and let $\mathcal{H}$ be the set of horoballs in $\widehat{\mathcal{H}}$ based at fixed points of parabolic elements of $\rho(\pi_1(M))$. Note that $\mathcal{H}$ is a precisely invariant system of horoballs for $\rho(\pi_1(M))$. Let $N^0 = (\mathbf{H}^3 - \mathcal{H})/\rho(\pi_1(M))$, $\widehat{N}^0 = (\mathbf{H}^3 - \widehat{\mathcal{H}})/\widehat{\Gamma}$, and $N_i = \mathbf{H}^3/\rho_i(\pi_1(M))$. Let $\{f_i\colon Z_i \to N_i\}$ be the sequence of orientation-preserving biLipschitz diffeomorphisms produced by Lemma 3.5.

*The strongly convergent case.* We first suppose that $\widehat{\Gamma} = \rho(\pi_1(M))$. Let $M_0$ be a compact core for $N$. Since $M_0 \subset Z_i$ for all large enough $i$, we may define $X_i = f_i(M_0)$. Lemma 3.6 implies that if $g_i$ is the restriction of $f_i$ to $M_0$, then $(g_i)_* = \rho_i \circ \rho^{-1}$, as maps of $\pi_1(M_0) = \rho(\pi_1(M))$ into $\rho_i(\pi_1(M))$, for all large enough $i$, and so $X_i$ is a compact core for $N_i$ for all large enough $i$. Therefore, if $\psi\colon M \to M_0$ is a homotopy equivalence such that $\psi_*$ is conjugate to $\rho$, then $\Psi(\rho_i) = [(X_i, g_i \circ \psi)]$ is equal to $\Psi(\rho) = [(M_0, \psi)]$ for all large enough $i$. Hence, Theorem A holds whenever $\widehat{\Gamma} = \rho(\pi_1(M))$.

Suppose now that $\Omega(\rho(\pi_1(M)))$ is empty. Since $N$ is topologically tame, Theorem 9.2 of Canary [11] guarantees that $\widehat{\Gamma} = \rho(\pi_1(M))$, and so again Theorem A holds.

For the remainder of the proof, we may assume that $\widehat{\Gamma}$ is not equal to $\rho(\pi_1(M))$ and that $\Omega(\rho(\pi_1(M))$ is nonempty.

*Decomposing the compact core.* Let $(M_0, P_0)$ be a relative compact core for $N^0$. Since $M_0$ is homotopy equivalent to $M$, $M_0$ also has nonempty incompressible boundary (using Theorem 7.1 of Hempel [15]). Moreover, since $P_0$ is a collection of incompressible annuli and tori in $\partial M_0$, each component of $\partial M_0 - P_0$ is an incompressible surface.

There exists a maximal collection $\mathcal{A} = \{A_1, \ldots, A_l\}$ of disjoint, nonparallel, essential annuli in $M_0$ such that, for each $i$, one boundary component of $A_i$ lies in $P_0$ and the other lies in a component of $\partial M_0 - P_0$. Select disjoint closed regular neighborhoods of the $A_i$, each intersecting $P_0$ in an annulus in the interior of $P_0$, and let $\mathcal{U} = \{\mathcal{U}(A_1), \ldots, \mathcal{U}(A_l)\}$ be the collection of their interiors. Set $M' = M_0 - \left(\cup_{i=1}^l \mathcal{U}(A_i)\right)$.

Let $R$ be a component of $M'$ and set $B = \overline{\partial R - \partial M_0}$. Every component of $B$ is an essential annulus which is properly homotopic to some element of $\mathcal{A}$. Notice that this implies that $\pi_1(R)$ injects into $\pi_1(M)$. In the next two paragraphs we show that no component of $M'$ is either a solid torus or a thickened torus.

Suppose that $R$ is a solid torus. Note that, by Lemma 3.1, any component of $P_0 \cap R$, hence any component of $B$, is a primitive annulus. If $B$ is empty, then $R = M_0$, contradicting the fact that $\pi_1(M_0) \cong \pi_1(M)$ is nonabelian. If $B$ has only one component, then $B$ is properly homotopic to an annulus with both boundary components in $P_0$, and Lemma 3.1 would show that $B$ is not essential. If $B$ has two components, they must be properly homotopic, and so must form the frontier of a single $\mathcal{U}(A_i)$. This would imply that either $M_0$ is a thickened torus, contradicting the fact that $\pi_1(M)$ is nonabelian, or that $M_0$ is an $I$-bundle over the Klein bottle, contradicting Lemma 3.1. If $B$ has more than two components, then there exists a properly embedded incompressible annulus $A$ in $R$ joining two components of $R \cap P_0$ such that both components of $V - A$ contain components of $B$. By Lemma 3.1, $A$ is properly homotopic into $P_0$, showing that $A$ is the frontier of a solid torus $W$ in $M_0$ such that $W \cap \partial M_0 \subset P_0$. By the choice of $A$, there is a component of $B$ in $W$, which contradicts the fact that each component of $B$ is essential. Therefore, no component of $M'$ is a solid torus.

Suppose that $R$ is a thickened torus. Since $\pi_1(M_0)$ is nonabelian, $R \neq M_0$ so that $B$ is nonempty. Let $T$ be a component of $\partial R$ containing an annulus component $A$ of $B$ that meets a component $P_0^1$ of $P_0$. By Lemma 3.1, $T$ is homotopic into a torus component of $P_0$, necessarily $P_0^1$. Also by Lemma 3.1, the other boundary component of $A$ must lie in $\partial M_0 - P_0^1$, so that $A$ is nonseparating. Let $\widehat{M_0}$ be the infinite cyclic cover formed by splitting $M_0$ along $A$ and laying copies end-to-end. Then $T$ lifts to $\widehat{M_0}$ and the homotopy of $T$ into $P_0^1$ lifts to a homotopy carrying the lifted $T$ into an open annulus that covers $P_0^1$. Therefore $\pi_1(T)$ is conjugate into an infinite cyclic subgroup of $\pi_1(M_0)$, which is impossible since $T$ is incompressible.

*Finding relative compact carriers.* Let $\{R_1, \ldots, R_n\}$ be the components of $M'$. For each $1 \leq j \leq n$, let $\Gamma_j = \pi_1(R_j) \subset \rho(\pi_1(M))$. (Notice that $\Gamma_j$ is really only well-defined up to conjugacy.) Since no $R_j$ is homeomorphic to either a solid torus or a thickened torus, each $\Gamma_j$ is nonabelian. Let $\mathcal{H}_j$ be the horoballs in $\widehat{\mathcal{H}}$ centered at parabolic fixed points of $\Gamma_j$. Then $\mathcal{H}_j$ is a precisely invariant set of horoballs for $\Gamma_j$. Set $\widetilde{N}_j = \mathbf{H}^3/\Gamma_j$ and $\widetilde{N}_j^0 = (\mathbf{H}^3 - \mathcal{H}_j)/\Gamma_j$, and let $p_j \colon \widetilde{N}_j \to N$ be the obvious covering map. Note that $R_j$ lifts to a relative compact core $\widetilde{R}_j$ for $\widetilde{N}_j$, and that $R_j$ is a relative compact carrier for $\Gamma_j$ in $N$. Let $P_j = R_j \cap \partial N^0$.

Notice that each component of $\partial R_j - P_j$ is obtained from a component of $\partial M_0 - (\cup_i \mathcal{U}(A_i))$ by appending collar neighborhoods of a collection (possibly empty) of its boundary components. Since $\partial M_0$ is incompressible and each component of $\partial M_0 \cap (\cup_i \mathcal{U}(A_i))$ is an incompressible annulus, each component of $\partial R_j - P_j$ is incompressible. By construction, $R_j$ contains no essential annulus with one boundary component in $P_j$ and the other in $\partial R_j - P_j$, and so Lemma 3.2 implies that $\Gamma_j$ is either a generalized web group or a degenerate group without accidental parabolics.

Renumber so that $\{\Gamma_1, \ldots, \Gamma_m\}$ are generalized web groups and so that $\{\Gamma_{m+1}, \ldots, \Gamma_n\}$ are degenerate groups. By construction, $\{\Gamma_1, \ldots, \Gamma_n\}$ are the vertex groups of a Bass-Serre graph for $\Gamma$ which has infinite cyclic edge groups. It follows, since each $\Gamma_j$ is nonabelian, that the groups $\{\Gamma_1, \ldots, \Gamma_n\}$ are nonconjugate. Since $\{R_1, \ldots, R_m\}$ is a disjoint collection of relative compact carriers for $\{\Gamma_1, \ldots, \Gamma_m\}$ and each component of the closure of $N^0 - R_j$ is noncompact, Lemma 4.1 guarantees that $\{\Gamma_1, \ldots, \Gamma_m\}$ is a precisely embedded system of generalized web subgroups of $\rho(\pi_1(M))$.

Suppose that $1 \le j, k \le m$, so that $\Gamma_j$ and $\Gamma_k$ are generalized web groups. If $\gamma \in \widehat{\Gamma} - \rho(\pi_1(M))$, then Proposition 3.7 guarantees that $\gamma(\Lambda(\Gamma_j)) \cap \Lambda(\Gamma_k)$ contains at most one point. Since the boundary of every component of $\Omega(\Gamma_j)$ and of $\Omega(\Gamma_k)$ is a Jordan curve, this implies that $\gamma(\Lambda(\Gamma_j))$ is contained in the closure of a component of $\Omega(\Gamma_k)$. We conclude that $\{\Gamma_1, \ldots, \Gamma_m\}$ is a precisely embedded system of generalized web subgroups of $\widehat{\Gamma}$. Hence, Proposition 4.2 guarantees that there exists a disjoint collection $\{\widehat{Y}_1, \ldots, \widehat{Y}_m\}$ of relative compact carriers for $\{\Gamma_1, \ldots, \Gamma_m\}$ in $\widehat{N}^0$.

We have just constructed relative compact carriers for the generalized web groups. In the next two paragraphs we construct relative compact carriers for the degenerate groups. Suppose that $m+1 \le j \le n$, so that $\Gamma_j$ is a degenerate group without accidental parabolics. Bonahon's theorem [8] implies that $\widetilde{N}^0_j$ is homeomorphic to $F_j \times \mathbf{R}$ (for some compact surface $F_j$) and that $\widetilde{R}_j$ may be identified with $F_j \times [-1, 0]$. Let $\widetilde{U}_j$ be the component of $\widetilde{N}^0_j - \widetilde{R}_j$ which is a neighborhood of the geometrically infinite end of $\widetilde{N}^0_j$. We may assume that $\widetilde{U}_j$ is identified with $F_j \times (0, \infty)$. Thurston's covering theorem, see [11], guarantees that the covering map $p_j \colon \widetilde{N}_j \to N$ is finite-to-one on $\widetilde{U}_j$. If $p_j$ is not injective on $\widetilde{U}_j$, then $\widetilde{U}_j$ contains a component of $p_j^{-1}(\text{Fr}(R_j))$ which is a properly embedded compact incompressible surface, and hence isotopic to $F_j \times \{0\}$. Thus, the component $\widetilde{B}$ of the closure of $\widetilde{U}_j - p_j^{-1}(R_j)$ which is adjacent to $\widetilde{R}_j$ is compact. However, this would imply that $B = p_j(\widetilde{B})$ is a compact component of $N^0 - R_j$ which is impossible. Thus $p_j$ is injective on $\widetilde{U}_j$ and $U_j = p_j(\widetilde{U}_j)$ is a neighborhood of a geometrically infinite end $E_j$ of $N^0$.

Since $N$ is topologically tame and $E_j$ is geometrically infinite, Proposition 3.8 implies that there exists a neighborhood $U'_j$ of $E_j$ that embeds in $\widehat{N}$. Since we may choose $U'_j$ to be identified with $F_j \times (k, \infty)$ (for some $k > 0$), we may find a relative compact carrier $\widehat{Y}_j$ for $\Gamma_j$ in $\widehat{N}$ which is disjoint from all the previously constructed relative compact carriers. Hence, we may successively extend $\{\widehat{Y}_1, \ldots, \widehat{Y}_m\}$ to a disjoint collection $\{\widehat{Y}_1, \ldots, \widehat{Y}_n\}$ of relative compact carriers for $\{\Gamma_1, \ldots, \Gamma_n\}$ in $\widehat{N}$. Notice that we may lift $\{\widehat{Y}_1, \ldots, \widehat{Y}_n\}$ to a disjoint collection $\{Y_1, \ldots, Y_n\}$ of relative compact carriers for $\{\Gamma_1, \ldots, \Gamma_n\}$ in $N^0$.

*Organizing the cusps.* It will be useful to characterize how cuspidal regions in $N$ cover cuspidal regions in $\widehat{N}$. Given a parabolic fixed point $x$ of $\rho(\pi_1(M))$, let $H$ be the horoball in $\widehat{\mathcal{H}}$ based at $x$. Then, $Q = \partial H / \Gamma_x$ is a component of $\partial N^0$, where $\Gamma_x = \mathrm{st}_{\rho(\pi_1(M))}(x)$. We note, by the construction of $\widehat{\mathcal{H}}$, that $\pi(Q) = \partial H / \widehat{\Gamma}_x$, where $\widehat{\Gamma}_x = \mathrm{st}_{\widehat{\Gamma}}(x)$.

If $Q$ is a torus, then $\widehat{\Gamma}_x$ and $\Gamma_x$ both have rank 2, so that $\Gamma_x$ has finite index in $\widehat{\Gamma}_x$. In this case, Lemma 3.4 implies that $\Gamma_x = \widehat{\Gamma}_x$, and so $Q$ embeds in $\widehat{N}$.

If $Q$ is an annulus, then $\Gamma_x$ has rank 1. If $\widehat{\Gamma}_x$ also has rank 1, then again Lemma 3.4 implies that $\Gamma_x = \widehat{\Gamma}_x$, and $Q$ embeds in $\widehat{N}$. If $\widehat{\Gamma}_x$ has rank 2, then $\pi(Q)$ is a torus and Lemma 3.4 implies that some core curve for $Q$ embeds in $\pi(Q)$.

We observe that if $Q$ and $Q'$ are distinct components of $\partial N^0$, then $\pi(Q)$ and $\pi(Q')$ are distinct components of $\partial \widehat{N}^0$. If not, then since $\widehat{\mathcal{H}}$ is a precisely invariant system of horoballs for $\widehat{\Gamma}$, we have that $\pi(Q) = \pi(Q')$. Let $Q = \partial H / \Gamma_x$ and $Q' = \partial H' / \Gamma_{x'}$. Since $\pi(Q) = \pi(Q')$, there exists $\gamma \in \widehat{\Gamma} - \rho(\pi_1(M))$ such that $\gamma(x) = x'$. Proposition 3.7 implies that $\gamma \Gamma_x \gamma^{-1}$ and $\Gamma_{x'}$ do not together generate a rank two abelian subgroup of $\widehat{\Gamma}$. Therefore, $\gamma \Gamma_x \gamma^{-1} \cap \Gamma_{x'}$ is nonempty. Let $\rho(g)$ and $\rho(g')$ be elements of $\Gamma_x$ and $\Gamma_{x'}$, respectively, such that $\gamma \rho(g) \gamma^{-1} = \rho(g')$. Writing $\gamma = \lim \gamma_i$, we see that $\gamma_i \rho_i(g) \gamma_i^{-1} = \rho_i(g')$ for large enough $i$ (see Lemma 3.6 in Jørgensen and Marden [21]). If $\gamma_i = \rho_i(h_i)$, then $h_i g h_i^{-1} = g'$, so $\rho(h_i)(x) = x'$, contradicting our assumption that $Q$ and $Q'$ are distinct components of $\partial N^0$.

Enumerate the components of $\partial N^0$ as $Q_1, \ldots, Q_s$, $Q_{s+1}, \ldots, Q_t$, $Q_{t+1}, \ldots, Q_u$, so that if $1 \leq k \leq s$, then $Q_k$ is an infinite annulus which covers a torus in $\widehat{N}$, if $s+1 \leq k \leq t$, then $Q_k$ is an infinite annulus which embeds in $\widehat{N}$, and if $t+1 \leq k \leq u$, then $Q_k$ is a torus, which necessarily embeds in $\widehat{N}$. Renumber the first $s$ as $Q_1, \ldots, Q_r, Q_{r+1}, \ldots, Q_s$ so that if $1 \leq k \leq r$, then $(\bigcup_j Y_j) \cap Q_k$ has at least three components, and if $r+1 \leq k \leq s$, then $(\bigcup_j Y_j) \cap Q_k$ has at most two components.

*Constructing collar neighborhoods.* The following explicit construction of a collar neighborhood of a submanifold of the boundary of a cuspidal region will be used in many of the constructions in the proof. Let $X = \mathbf{H}^3/\Theta$ be a hyperbolic 3-manifold and let $\mathcal{H}_X$ be a precisely invariant system of horoballs for $\Theta$. If $C$ is a component of $\mathcal{H}_X/\Theta$, then there exists a homeomorphism $\widehat{r} \colon C \to \partial C \times [0, \infty)$ given by $\widehat{r}(x) = (y, t)$, where $y$ is the orthogonal projection of $x$ onto $\partial C$ and $t$ is the hyperbolic distance between $x$ and $y$. If $Z$ is a submanifold of $\partial C$, set $\mathcal{N}(Z) = \widehat{r}^{-1}(Z \times [0, 1])$.

*A new compact core for the algebraic limit.* Let $M_0^+ = M_0 \cup \mathcal{N}(P_0)$. Since $M_0^+$ is obtained from $M_0$ by appending collar neighborhoods of each component of $P_0$, $M_0^+$ is also a compact core for $N$. Let $M_0' = M' \cup \mathcal{N}(P_0)$. Then $M_0^+$ is obtained from $M_0'$ by appending $\mathcal{U}$. Each component $U$ of $\mathcal{U}$ is a solid torus attached to $M_0'$ along an annulus which is primitive in $U$. Thus, since $M_0^+$ is a compact core for $N$, so is $M_0'$.

If $1 \le k \le t$, so that $Q_k$ is an annulus, let $C_k$ be the minimal annulus in $Q_k$ which contains $(\cup_j Y_j) \cap Q_k$. If $t+1 \le k \le u$, so that $Q_k$ is a torus, let $C_k = Q_k$.

Let $M_1 = (\cup_{j=1}^n Y_j) \cup (\cup_{k=1}^u \mathcal{N}(C_k))$. We claim that $M_1$ is also a compact core for $N$. We establish our claim by constructing a homotopy equivalence $h \colon M_0' \to M_1$ and a homotopy $L \colon M_0' \times I \to N$ from the inclusion $i_{M_0'}$ to $i_{M_1} \circ h$.

Lemma 3.3 implies that there are homotopies $H_j \colon R_j \times I \to N^0$ such that $(H_j)_0$ is the inclusion, $(H_j)_1$ carries $R_j$ homeomorphically to $Y_j$, and $H_j((R_j \cap \partial N^0) \times I) \subset \partial N^0$. Define a partial homotopy

$$J \colon P_0 \times \{0\} \cup ((\cup_{j=1}^n R_j) \cap P_0) \times I \to \partial N^0$$

by letting $J_0$ equal the inclusion and using the restrictions of the homotopies $H_j$ on $((\cup_{j=1}^n R_j) \cap P_0) \times I$. By the homotopy extension property, $J$ extends to all of $P_0 \times I$. This extends, using the $H_j$, to a homotopy $K \colon (P_0 \cup (\cup_{j=1}^n R_j)) \times I \to N^0$ which carries $P_0 \times I$ into $\partial N^0$.

Let $G \colon \partial N^0 \times I \to \partial N^0$ be a deformation retraction onto $\cup_{k=1}^u C_k$, i.e. $G_0$ is the identity, $G_1$ is a retraction onto $\cup_{k=1}^u C_k$, and each $G_t$ restricts to the identity on $\cup_{k=1}^u C_k$. Extend $G$ using the identity maps on $\cup_{j=1}^n Y_j$ to a deformation retraction $G \colon (\partial N^0 \cup (\cup_{j=1}^n Y_j)) \times I \to \partial N^0 \cup (\cup_{j=1}^n Y_j)$ onto $(\cup_{k=1}^u C_k) \cup (\cup_{j=1}^n Y_j)$. By the homotopy extension property, $G$ extends to a homotopy $G \colon N^0 \times I \to N^0$ with $G_0$ equal to the identity map.

The homotopy $L \colon \left(P_0 \cup (\cup_{j=1}^n R_j)\right) \times I \to N^0$ given by letting $L_t = G_t \circ K_t$ for all $t$, starts at the inclusion of $P_0 \cup (\cup_{j=1}^n R_j)$ into $N^0$ and ends at a map carrying $P_0 \cup (\cup_{j=1}^n R_j)$ into $(\cup_{k=1}^u C_k) \cup (\cup_{j=1}^n Y_j)$. The map $L_1$ is a homeomorphism on each $R_j$ (since it agrees with $(H_j)_1$) and carries $P_0$ into $\cup_{k=1}^u C_k$. Since $L_t(\partial N^0) \subseteq \partial N^0$ for each $t$, we may extend $L$ over $\mathcal{N}(P_0) \times I$ by setting $L(\widehat{r}^{-1}(x,s),t) = \widehat{r}^{-1}(L(x,t),s)$. Define $h = i_{M_1}^{-1} \circ L_1$. Lemma 5.2 implies that $h$ is a homotopy equivalence from $M_0'$ to $M_1$, so that $L$ is the desired homotopy.

Since $M_0'$ is a compact core and $i_{M_1} \circ h$ is homotopic to $i_{M_0'}$, $M_1$ is also a compact core. Hence, $\Psi(\rho) = [(M_1, \psi)]$, where $\psi \colon M \to M_1$ is a homotopy equivalence such that $\psi_* = \rho$.

*A compact submanifold of the geometric limit.* We similarly construct a compact submanifold $M_2$ of $\widehat{N}$ which is homotopy equivalent to $M_1$ and which will pull back (via the biLipschiz diffeomorphism $f_i$) to give a compact core for $N_i$, for all large enough $i$. If $1 \leq k \leq s$, so that $Q_k$ is an annulus and $\pi(Q_k)$ is a torus, let $B_k$ be any annulus in $\pi(Q_k)$ which contains $(\bigcup_j \widehat{Y}_j) \cap \pi(Q_k)$ and is contained within $\pi(C_k)$. If $s+1 \leq k \leq u$, so that $Q_k$ embeds under $\pi$, let $B_k = \pi(C_k)$. Set $M_2 = (\bigcup_j \widehat{Y}_j) \cup (\bigcup_k \mathcal{N}(B_k))$.

*The primitive shuffle equivalence.* We now define an explicit homotopy equivalence, $\phi\colon M_1 \to M_2$, which will turn out to be our desired primitive shuffle. Define $\phi$ to agree with the covering map $\pi$ on each $Y_j$. If $s+1 \leq k \leq u$, so that $Q_k$ embeds under $\pi$, define $\phi$ to agree with $\pi$ on $\mathcal{N}(C_k)$. If $r+1 \leq k \leq s$, choose the restriction of $\phi$ to $\mathcal{N}(C_k)$ to be an orientation-preserving homeomorphism of $\mathcal{N}(C_k)$ to $\mathcal{N}(B_k)$ which agrees with $\pi$ on every component of $(\cup_j Y_j) \cap C_k$. If $1 \leq k \leq r$, we may only choose $\phi$ to be a homotopy equivalence of $\mathcal{N}(C_k)$ to $\mathcal{N}(B_k)$ which agrees with $\pi$ on $(\cup_j Y_j) \cap C_k$. Lemma 5.2 implies that $\phi$ is a homotopy equivalence.

It remains to show that $\phi$ is a primitive shuffle equivalence. If $1 \leq k \leq r$, then the submanifold $\mathcal{N}(C_k)$ is homeomorphic to a solid torus and admits a Seifert fibering so that $\mathcal{N}(C_k) \cap \partial M_1$ is a collection of at least three fibered annuli. Since each component of the frontier of $\mathcal{N}(C_k)$ is an essential annulus in $M_1$, the inclusion map is an admissible, essential embedding. Hence, $\mathcal{N}(C_k)$ is admissibly isotopic into a component $V_k$ of $\Sigma(M_1)$. In fact, we may assume that $\mathcal{N}(C_k) \subset V_k$. Since there are essential annuli in $V_k$ which are homotopic but not parallel, $V_k$ cannot be an $I$-bundle component of $\Sigma(M_1)$. Since $\pi_1(C_k)$ is a maximal abelian subgroup of $\pi_1(M_1)$, $V_k$ cannot be homeomorphic to a thickened torus. Moreover, again since $\pi_1(C_k)$ is a maximal abelian subgroup of $\pi_1(M_1)$, $V_k$ is a primitive solid torus. Proposition 6.1 then implies, perhaps after isotoping $\Sigma(M_2)$, that $W_k = \phi(V_k)$ is a primitive solid torus component of $\Sigma(M_2)$ for all $1 \leq k \leq r$. By construction, $\phi$ is an orientation-preserving homeomorphism from the closure of $M_1 - \cup_{k=1}^r V_k$ to the closure of $M_2 - \cup_{k=1}^r W_k$. Therefore, $\phi$ is a primitive shuffle.

*Compact cores for the approximates.* We now show that $X_i = f_i(M_2)$ is a compact core for $N_i$ for all large enough $i$. We begin by studying the images of the tori in $\widehat{N}$ which arise as the projections of boundaries of rank one cusps of $N$. If $1 \leq k \leq s$, let $T_k = \pi(Q_k)$. Choose a meridian-longitude system $(a_k, b_k)$ for $\pi_1(T_k)$ so that the longitude $b_k$ is a core curve for $B_k$ and the meridian $a_k$ intersects $b_k$ exactly once. We recall that $\{f_i\colon Z_i \to N_i\}$ is the sequence of orientation-preserving biLipschitz diffeomorphisms produced by Lemma 3.5. For all large enough $i$, $M_2$, $T_1, \ldots, T_{s-1}$ and $T_s$ all lie in $Z_i$. Let $X_i = f_i(M_2)$, and, for $1 \leq k \leq s$, let $T_k^i = f_i(T_k)$. Assign the meridian-longitude system $(a_k^i = f_i(a_k), b_k^i = f_i(b_k))$ to $T_k^i$.

Let $\widetilde{b}_k$ be a lift of $b_k$ to $C_k$. Lemma 3.6 implies that, for all large enough $i$, $b_k^i = \rho_i(\rho^{-1}(\widetilde{b}_k))$ as elements of $\rho_i(\pi_1(M))$. Since $\widetilde{b}_k$ generates a maximal abelian subgroup of $\Gamma$, $b_k^i$ generates a maximal abelian subgroup of $\rho_i(\pi_1(M))$. Therefore, $T_k^i$ bounds a solid torus $V_k^i$ for all large enough $i$. For the remainder of the argument assume that $i$ has been chosen large enough that each $T_k^i$ bounds a solid torus and that $b_k^i$ generates a maximal abelian subgroup of $\rho_i(\pi_1(M))$.

Let $X_i^+ = X_i \cup (\cup_{k=1}^s V_k^i)$ and let $X_i^- = f_i(M_2 \cap \widehat{N}^0)$. Lemma 3.6 implies that, for all large enough $i$, $(f_i \circ \pi)_*$ agrees with $\rho_i \circ \rho^{-1}$ on $\pi_1(Y_j)$ for all $j$. Since $\pi$ embeds $Y_j$ into $M_2 \cap \widehat{N}^0$ and $\pi_1(Y_j)$ is nonabelian, we may conclude that $\pi_1(X_i^-)$ is nonabelian for all large enough $i$. Therefore, $V_k^i$ lies on the opposite side of $T_k^i$ from $X_i^-$ for all large enough $i$. Moreover, if $A$ is any component of the frontier of $V_k^i$ in $X_i^+$, then $A$ is an incompressible annulus in $f_i(B_k)$. Therefore, the image of the inclusion of $\pi_1(A)$ into $\rho_i(\pi_1(M))$ is generated by $b_k^i$. Since $b_k^i$ generates a maximal abelian subgroup of $\rho_i(\pi_1(M))$, we see that each $V_k^i$ is a primitive solid torus in $X_i^+$. Moreover, since each $V_k^i \cap X_i = f_i(\mathcal{N}(B_k))$ is a primitive solid torus containing the frontier of $V_k^i$ in $X_i^+$, $X_i$ is a compact core for $X_i^+$.

Since $\pi(\mathcal{N}(C_k))$ lies on the opposite side of $T_k$ from $M_2 \cap \widehat{N}^0$, $f_i(\pi(\mathcal{N}(C_k)))$ must lie entirely in $V_k^i$ for all $1 \le k \le s$ and all large enough $i$. Thus, $f_i(\pi(M_1))$ is contained entirely in $X_i^+$ for all large enough $i$. If $\pi' : M_1 \to \pi(M_1)$ is the restriction of $\pi$ to $M_1$, then Lemma 3.6 implies that $(f_i \circ \pi')_* = \rho_i \circ \rho^{-1}$, as maps of $\pi_1(M_1) = \rho(\pi_1(M))$ into $\rho_i(\pi_1(M))$, for all large enough $i$.

Since $\cup_{k=1}^s \mathcal{N}(C_k)$ and $\cup_{i=1}^s V_k^i$ are both collections of primitive solid tori, $f_i \circ \pi'$ is a homeomorphism from the closure of $M_1 - \cup_{k=1}^s \mathcal{N}(C_k)$ to the closure of $X_i^+ - \cup_{k=1}^s V_k^i$, and $f_i \circ \pi'$ takes $\cup_{k=1}^s \mathcal{N}(C_k)$ into $\cup_{i=1}^s V_k^i$, Lemma 5.2 implies that $f_i \circ \pi'$ is a homotopy equivalence of $M_1$ to $X_i^+$ for all large enough $i$. By construction, $f_i \circ \phi$ is a homotopy equivalence from $M_1$ to $X_i$. As $X_i$ is a compact core for $X_i^+$, $f_i \circ \phi$ is also a homotopy equivalence of $M_1$ to $X_i^+$. Since $f_i \circ \phi$ and $f_i \circ \pi'$ both restrict to the same homeomorphism from the closure of $M_1 - \cup_{k=1}^s \mathcal{N}(C_k)$ to the closure of $X_i^+ - \cup_{k=1}^s V_k^i$ (and each takes $\cup_{k=1}^s \mathcal{N}(C_k)$ into $\cup_{k=1}^s V_k^i$), Lemma 7.1 guarantees that there exists an orientation-preserving homeomorphism $r_i : M_1 \to M_1$ (which is the identity off of $\cup_{k=1}^s \mathcal{N}(C_k)$) such that $f_i \circ \phi \circ r_i$ and $f_i \circ \pi'$ are homotopic as maps of $M_1$ into $X_i^+$. Therefore, for all large enough $i$, $f_i \circ \phi \circ r_i$ is a homotopy equivalence from $M_1$ to $X_i$ with the property that $(f_i \circ \phi \circ r_i)_* = \rho_i \circ \rho^{-1}$ (as maps of $\pi_1(M_1) = \rho(\pi_1(M))$ into $\pi_1(N_i) = \rho_i(\pi_1(M))$). It follows that $X_i$ is a compact core for $N_i$ for all large enough $i$.

*The marked homeomorphism types are primitive shuffle equivalent.* Recall that we have written $\Psi(\rho) = [(M_1, \psi)]$ and that we have described, for all large enough $i$, a homotopy equivalence $f_i \circ \phi \circ r_i$ from $M_1$ to a compact core $X_i$ for $N_i$ such that $(f_i \circ \phi \circ r_i)_* = \rho_i \circ \rho^{-1}$. Hence, again for all large enough $i$,

$$\Psi(\rho_i) = [(X_i, f_i \circ \phi \circ r_i \circ \psi)] = [(M_2, \phi \circ r_i \circ \psi)].$$

(The second equality follows from the fact that $f_i$ restricts to an orientation-preserving homeomorphism of $M_2$ to $X_i$.) Since $\phi \circ r_i$ and $\phi$ both restrict to the same homeomorphism of the closure of $M_1 - \cup_{k=1}^{s} \mathcal{N}(C_k)$ to the closure of $M_2 - \cup_{k=1}^{s} \mathcal{N}(B_k)$, Lemma 7.1 implies that there exists an orientation-preserving homeomorphism $r_i'\colon M_2 \to M_2$ so that $\phi \circ r_i$ is homotopic to $r_i' \circ \phi$. Hence,

$$\Psi(\rho_i) = [(M_2, \phi \circ r_i \circ \psi)] = [(M_2, r_i' \circ \phi \circ \psi)] = [(M_2, \phi \circ \psi)]$$

for all large enough $i$. Therefore, $\Psi(\rho_i)$ is primitive shuffle equivalent to $\Psi(\rho)$ for all large enough $i$. This completes the proof of our initial claim and thus establishes Theorem A. □

In establishing the corollaries it will be useful to consider the map $\widehat{\Psi}\colon AH(\pi_1(M)) \to \widehat{\mathcal{A}}(M)$ which is obtained by composing $\Psi\colon AH(\pi_1(M)) \to \mathcal{A}(M)$ with the quotient map $q\colon \mathcal{A}(M) \to \widehat{\mathcal{A}}(M)$. In this language, Theorem A asserts that $\widehat{\Psi}$ is continuous, where we give $\widehat{\mathcal{A}}(M)$ the discrete topology.

COROLLARY 8.1. *Let $M$ be a compact, hyperbolizable 3-manifold with nonempty incompressible boundary. Then the map $\widehat{\Psi}\colon AH(\pi_1(M)) \to \widehat{\mathcal{A}}(M)$ is continuous (i.e. locally constant).*

The proof of Theorem A has the following corollary which gives further geometric restrictions on when the marked homeomorphism type can change in the algebraic limit. In particular, we see that if the marked homeomorphism type of the algebraic limit differs from that of the approximates, then there must be a specific type of "new" rank two abelian subgroup of the geometric limit.

COROLLARY 8.2. *Let $M$ be a compact, hyperbolizable 3-manifold with nonempty incompressible boundary. Suppose that the sequence $\{\rho_i\}$ in $\mathcal{D}(\pi_1(M))$ converges to $\rho$ and $\{\rho_i(\pi_1(M))\}$ converges geometrically to $\widehat{\Gamma}$. If there does not exist a primitive solid torus component $V$ of $\Sigma(M)$ such that $V \cap \partial M$ has at least three components and $\rho(\pi_1(V))$ lies in a rank-two free abelian subgroup of $\widehat{\Gamma}$, then $\Psi(\rho_i) = \Psi(\rho)$ for all sufficiently large $i$.*

## 9. Lifting primitive shuffles

In the proof of Theorem B, we will need to know that any primitive shuffle is homotopic to one that misses the core curves of the primitive solid tori in the image, and which lifts to an embedding in some covering space of the complement of the core curves.

PROPOSITION 9.1. *Let $s: M_1 \to M_2$ be a primitive shuffle with respect to $\mathcal{V}_1$ and $\mathcal{V}_2$, and fix a collection $\Delta_2$ of core curves for the solid tori of $\mathcal{V}_2$. Then there exist a homotopy equivalence $r: M_1 \to M_2$ homotopic to $s$ and a covering space $\widetilde{M}_2$ of $M_2 - \Delta_2$ such that $r(M_1) \subseteq M_2 - \Delta_2$ and $r$ lifts to an embedding $\widetilde{r}: M_1 \to \widetilde{M}_2$.*

*Proof of Proposition* 9.1. Assume first that $\mathcal{V}_1$ consists of a single solid torus $V$, then $\mathcal{V}_2$ is also a single solid torus $W$. Let $S_1, \ldots, S_m$ be the components of $\overline{M_2 - W}$. Let $B_1, \ldots, B_n$ be the annuli of the frontier of $W$, with notation chosen so that the $B_i$ are cyclically ordered as one travels around an oriented meridian curve $\gamma_2$ in $\partial W$ that intersects each $B_i$ in a single arc.

The restriction of $s$ is an orientation-preserving homeomorphism $h$ from $\overline{M_1 - V}$ to $\overline{M_2 - W}$. For each $i$, let $A_i = h^{-1}(B_i)$, a component of the frontier of $V$. Orient the arcs $h^{-1}(\gamma_2 \cap B_i)$ by transferring the orientation of $\gamma_2$ using $h^{-1}$. Since $h$ preserves orientation, the arcs $h^{-1}(\gamma_2 \cap B_i)$ lie in a single oriented meridian curve $\gamma_1$ of $\partial V$, in such a way that their orientations agree with the orientation of $\gamma_1$. To see this, fix a generator of $\pi_1(V)$ and regard its image under $s_*$ as the distinguished generator of $\pi_1(W)$. Let $\beta_i$ be the boundary component of $B_i$ that contains the initial point of $\gamma_2 \cap B_i$, then $h^{-1}(\beta_i)$ is the boundary component of $A_i$ that contains the initial point of $h^{-1}(\gamma_2 \cap B_i)$. Orient the $\beta_i$ to represent the distinguished generator of $\pi_1(W)$; under $h^{-1}$ these orientations transfer to orientations of $h^{-1}(\beta_i)$ which represent the distinguished generator of $\pi_1(V)$. For a given $B_i$, let $S_k$ be the component of $\overline{M_2 - W}$ that contains $B_i$. The pair consisting of the direction of $\gamma_2 \cap B_i$ followed by the direction of $\beta_i$ determines an orientation for the component of $\partial S_k$ that contains $B_i$ and hence an orientation for $S_k$. Since for all $i$ these pairs determine the same orientation of $\partial W$, the resulting orientations of the $S_k$ must either all agree or all disagree with the orientation of $M_2$. If they disagree, take the other choice of distinguished generator of $\pi_1(W)$ (and hence of $\pi_1(V)$) to ensure that they all agree. Since $h$ is orientation-preserving, the orientations on the components of $\mathrm{Fr}(V)$ determined by the pair consisting of the direction of $h^{-1}(\gamma_2 \cap B_i)$ followed by the direction of $h^{-1}(\beta_i)$ all agree with the orientation on $\mathrm{Fr}(V)$ determined by the orientation on $\overline{M_1 - V}$. Consequently, they all determine the same orientation of $\partial V$. Since the oriented $h^{-1}(\beta_i)$ all represent the distinguished generator of $\pi_1(V)$, this shows that if the endpoints of the $h^{-1}(\gamma_2 \cap B_i)$ are connected by arcs in $\overline{\partial V - \cup A_i}$, one obtains a meridian curve $\gamma_1$ which has an orientation agreeing with that of each $h^{-1}(\gamma_2 \cap B_i)$.

Define a permutation $\sigma$ of $\{1, \ldots, n\}$ by the rule that $\sigma(1) = 1$ and that $\gamma_1$ passes in order through $A_{\sigma(1)}, A_{\sigma(2)}, \ldots, A_{\sigma(n)}$.

Now we build the covering space $\widetilde{M}_2$ of $M_2 - \Delta_2$. Let $p_0: \widetilde{W} \to W - \Delta_2$ be the infinite cyclic covering of $W - \Delta_2$ such that the $B_i$ lift homeomorphically to $\widetilde{W}$. The lifts of the $B_i$ form an ordered infinite sequence

$$\ldots, B_n^{-1}, B_1^0, B_2^0, \ldots, B_n^0, B_1^1, \ldots, B_n^1, B_1^2, \ldots$$

of annuli as one travels along the lift $\widetilde{\gamma}_2$ of $\gamma_2$, where the $B_i^j$ map to $B_i$. For $1 \leq i \leq n$, put $C_{\sigma(i)} = B_{\sigma(i)}^i$. Thus $C_i$ projects to $B_i$, and $C_{\sigma(1)}, C_{\sigma(2)}, \ldots, C_{\sigma(n)}$ occur in order as one travels along $\widetilde{\gamma}_2$.

Attach copies $\widetilde{S}_k$ of the $S_k$ to $\widetilde{W}$ along $\cup C_i$ in such a way that $p_0$ extends, using orientation-preserving homeomorphisms from $\widetilde{S}_k$ to $S_k$, to a locally injective map $p_1 \colon \widetilde{W} \cup (\cup \widetilde{S}_k) \to M_2$. Note that $(p_1|_{\cup \widetilde{S}_k})^{-1} \circ s|_{\overline{M_1 - V}}$ is the unique lifting of $s|_{\overline{M_1 - V}}$ to an embedding into $\widetilde{W} \cup (\cup \widetilde{S}_k)$, and that it sends $A_i$ to $C_i$.

Choose an embedded solid torus $T$ in $\widetilde{W}$ that meets $\partial \widetilde{W}$ in $\cup C_i$. Since the $C_i$ occur in the order $C_{\sigma(1)}, C_{\sigma(2)}, \ldots, C_{\sigma(n)}$ on $\partial \widetilde{W}$, they occur in the same order as one travels around $\partial T$. Since the $A_{\sigma(1)}, \ldots, A_{\sigma(n)}$ occur in order as one travels around $V$, one can use a homeomorphism sending $V$ onto $T$, to extend $(p_1|_{\cup \widetilde{S}_k})^{-1} \circ s|_{\overline{M_1 - V}}$ to an embedding $\widetilde{r}_0 \colon M_1 \to \widetilde{W} \cup (\cup \widetilde{S}_k)$.

We will now enlarge $\widetilde{W} \cup (\cup \widetilde{S}_k)$ to a covering space $\widetilde{M}_2$ of $M_2 - \Delta_2$. Fix a $B_i^j$ which is not equal to $C_i$. Since $B_i$ is incompressible, $\pi_1(B_i)$ is a subgroup of $\pi_1(M_2 - \Delta_2)$. Let $q_i' \colon M_2' \to M_2 - \Delta_2$ be the covering space corresponding to this subgroup. There is an annulus $B_i' \subset M_2'$, mapped homeomorphically by $q_i'$ to $B_i$, such that $\pi_1(B_i') \to \pi_1(M_2')$ is an isomorphism. In particular, $B_i'$ separates $M_2'$. For one of the components of the complement, points near $B_i'$ map into $\overline{M_2 - W}$; denote by $M_2''$ the closure of this component, and let $q_i'' \colon M_2'' \to M_2 - \Delta_2$ be the restriction of $q_i'$. Attach $M_2''$ to $\widetilde{W} \cup (\cup \widetilde{S}_k)$ by identifying $B_i'$ to $B_i^j$ using the homeomorphism $\left(p_1|_{B_i^j}\right)^{-1} \circ \left(q_i''|_{B_i'}\right)$. Then $p_1$ extends to a local homeomorphism on the union by use of $q_i''$ on $M_2''$. Repeating this process for every $B_i^j$ that is not a $C_i$, we obtain $p \colon \widetilde{M}_2 \to M_2 - \Delta_2$. One sees easily that points in $M_2 - \Delta_2$ have evenly covered neighborhoods, so $p$ is a covering map.

Since $\widetilde{r}_0$ is an embedding into $\widetilde{W} \cup (\cup \widetilde{S}_k)$, we may regard it as an embedding of $M_1$ into $\widetilde{M}_2$. Let $r_0 \colon M_1 \to M_2$ be the composition $p \circ \widetilde{r}_0$. It carries $M_1$ into $M_2 - \Delta_2$ and lifts to the embedding $\widetilde{r}_0$ into $\widetilde{M}_2$, but it need not be homotopic to $s$. Notice that $(r_0)^{-1}(W) = V$ and $r_0|_{\overline{M_1 - V}}$ is the homeomorphism $s|_{\overline{M_1 - V}}$.

Applying Lemma 7.1, we can change $r_0$ by precomposing by a product of Dehn twists about frontier annuli of $V$ to obtain a map $r$ homotopic to $s$ as maps into $M_2$. Since $r$ still lifts to an embedding into $\widetilde{M}_2$, this completes the proof in the case where $\mathcal{V}_1$ and $\mathcal{V}_2$ have one component.

Now suppose that $\mathcal{V}_1$ has more than one component. Let $V_1, \ldots, V_n$ be the components of $\mathcal{V}_1$, and let $W_i$ be the component of $\mathcal{V}_2$ that contains $s(V_i)$. Start with infinite cyclic coverings $\widetilde{W}_i$ of the $W_i - \Delta_i$, and using $V_i$ select sequences $C_{i,1}, \ldots, C_{i,n_i}$ in the $\partial \widetilde{W}_i$ exactly as before. Again attach pieces $\widetilde{S}_k$ to $\cup \widetilde{W}_i$ along the $C_{i,j}$, and select solid tori $T_i$ in $\widetilde{W}_i$ meeting $\partial \widetilde{W}_i$ in $\cup_j C_{i,j}$. The construction of $\widetilde{M}_2$ and $\widetilde{r}_0$ proceeds as before, and again Lemma 7.1 provides a correction to produce a map $r$ which is homotopic to $s$. □

## 10. Dehn filling

In the proof of Theorem B, we make extensive use of Comar's generalization of Thurston's hyperbolic Dehn surgery theorem to the setting of geometrically finite Kleinian groups (see also Bonahon-Otal [9] and Bromberg [10]). We will apply this result to manifolds obtained from a hyperbolizable 3-manifold with incompressible boundary by removing core curves of primitive solid torus components of its characteristic submanifold. In this section, we discuss the operation of Dehn filling with emphasis on these special cases.

Let $\widehat{M}$ be a compact, hyperbolizable 3-manifold, and let $T_1, \ldots, T_k$ be a collection of toroidal boundary components of $\widehat{M}$. On each $T_i$, choose a meridian-longitude system $(m_i, l_i)$. Given a solid torus $V$ and a pair $(p_i, q_i)$ of relatively prime integers, we may form a new manifold by attaching $V$ to $\widehat{M}$ by an orientation-reversing homeomorphism $g\colon \partial V \to T_i$ so that, if $c$ is the meridian of $V$, then $g(c)$ is a $(p_i, q_i)$ curve on $T_i$ with respect to the chosen meridian-longitude system. We say that this manifold is obtained from $\widehat{M}$ by $(p_i, q_i)$-*Dehn filling along* $T_i$. Given a $k$-tuple $(\mathbf{p}, \mathbf{q}) = ((p_1, q_1), \ldots, (p_k, q_k))$ of pairs of relatively prime integers, let $\widehat{M}(\mathbf{p}, \mathbf{q})$ denote the manifold obtained by doing $(p_i, q_i)$-Dehn filling along $T_i$ for each $i$.

THEOREM 10.1 (Comar [14]). *Let $\widehat{M}$ be a compact, hyperbolizable 3-manifold and let $T = \{T_1, \ldots, T_k\}$ be a nonempty collection of tori in the boundary of $\widehat{M}$. Let $\widehat{N} = \mathbf{H}^3/\widehat{\Gamma}$ be a geometrically finite hyperbolic 3-manifold and let $\psi\colon \mathrm{int}(\widehat{M}) \to \widehat{N}$ be an orientation-preserving homeomorphism. Further assume that every parabolic element of $\widehat{\Gamma}$ lies in a rank-two parabolic subgroup. Let $(m_i, l_i)$ be a meridian-longitude basis for $T_i$. Let $\{(\mathbf{p}_n, \mathbf{q}_n) = ((p_n^1, q_n^1), \ldots, (p_n^k, q_n^k))\}$ be a sequence of $k$-tuples of pairs of relatively prime integers such that, for each $j$, $\{(p_n^j, q_n^j)\}$ converges to $\infty$ as $n \to \infty$.*

*For all sufficiently large $n$, there exists a representation $\beta_n\colon \widehat{\Gamma} \to \mathrm{PSL}_2(\mathbf{C})$ with discrete image such that*

1. *$\beta_n(\widehat{\Gamma})$ is geometrically finite, uniformizes $\widehat{M}(\mathbf{p}_n, \mathbf{q}_n)$, and every parabolic element of $\beta_n(\widehat{\Gamma})$ lies in a rank-two parabolic subgroup,*

2. *the kernel of $\beta_n \circ \psi_*$ is normally generated by $\{m_1^{p_n^1} l_1^{q_n^1}, \ldots, m_k^{p_n^k} l_k^{q_n^k}\}$, and*

3. *$\{\beta_n\}$ converges to the identity representation of $\widehat{\Gamma}$.*

*Moreover, if $i_n \colon \widehat{M} \to \widehat{M}(\mathbf{p}_n, \mathbf{q}_n)$ denotes the inclusion map, then for each $n$, there exists an orientation-preserving homeomorphism*

$$\psi_n \colon \mathrm{int}(\widehat{M}(\mathbf{p}_n, \mathbf{q}_n)) \to \mathbf{H}^3/\beta_n(\widehat{\Gamma})$$

*such that $\beta_n \circ \psi_*$ is conjugate to $(\psi_n)_* \circ (i_n)_*$.*

In order to apply the hyperbolic Dehn surgery theorem to manifolds obtained from a hyperbolizable manifold by removing core curves of solid torus components of its characteristic submanifold, we must first show that such manifolds are themselves hyperbolizable.

LEMMA 10.2. *Let $M$ be a compact, hyperbolizable 3-manifold with nonempty incompressible boundary, and let $\mathcal{V}$ be a collection of solid torus components of the characteristic submanifold $\Sigma(M)$ of $M$. Let $\Delta$ a collection of core curves of tori in $\mathcal{V}$, and suppose that $\mathcal{N}(\Delta)$ is an open regular neighborhood of $\Delta$ with $\overline{\mathcal{N}(\Delta)}$ in the interior of $\mathcal{V}$. Then, $\widehat{M} = M - \mathcal{N}(\Delta)$ is hyperbolizable. Moreover, $\Sigma(M) - \mathcal{N}(\Delta)$ is a characteristic submanifold for $\widehat{M}$.*

*Proof of Lemma* 10.2. Recall that Thurston's geometrization theorem (see, for example, Morgan [33] or Otal [36]) asserts that a compact 3-manifold with nonempty boundary is hyperbolizable if and only if it is irreducible and atoroidal.

Let $S$ be an embedded 2-sphere in $\widehat{M}$. Since $\widehat{M} \subset M$ and since $M$ is irreducible, $S$ bounds a 3-ball $B$ in $M$. As each component of $\mathcal{N}(\Delta)$ contains a closed homotopically nontrivial curve in $M$, $B$ cannot contain any component of $\mathcal{N}(\Delta)$. Therefore, $S$ bounds a ball in $\widehat{M}$, namely $B$, and so $\widehat{M}$ is irreducible.

We now show that $\Sigma(M) - \mathcal{N}(\Delta)$ is characteristic in $\widehat{M}$. This is an instance of a general principle discussed in Section 2.8 of [12]. This principle asserts that if one removes disjoint regular neighborhoods of finitely many fibers from a characteristic submanifold, the resulting submanifold is characteristic in the complement of the removed neighborhoods. For our specific case, however, since a brief and self-contained argument can be given using only basic properties of the characteristic submanifolds, we include it here.

First note that if an annulus in $\widehat{M}$ with boundary in $\partial M$ is essential in $M$, then it is essential in $\widehat{M}$. Consequently the frontier of $\Sigma(M) - \mathcal{N}(\Delta)$ is essential in $\widehat{M}$, so that $\Sigma(M) - \mathcal{N}(\Delta)$ is an essential fibered submanifold in $\widehat{M}$. Therefore it is isotopic into $\Sigma(\widehat{M})$, and we may assume that it is contained in the topological interior of $\Sigma(\widehat{M})$. Since $\mathrm{Fr}(\Sigma(M))$ consists of essential annuli, Corollary 5.7 of Johannson [19] shows that $\Sigma(\widehat{M})$ has an admissible fibering for which $\mathrm{Fr}(\Sigma(M))$ is a union of fibers. Consequently, each component $R$ of the closure of $\Sigma(\widehat{M}) - (\Sigma(M) - \mathcal{N}(\Delta))$ is an admissibly fibered submanifold of $M$.

We claim that each component $F$ of $\mathrm{Fr}(\Sigma(\widehat{M}))$ is an essential annulus in $M$. Note first that $F$ is incompressible in $M$, since any compressing disc for $F$ in $M$ could be surgered off the incompressible surface $\mathrm{Fr}(\Sigma(M))$ resulting in a compression in $\widehat{M}$. If $F$ is a torus, then since $M$ is atoroidal, $F$ is parallel in $M$ to a torus boundary component of $M$. Since the components of $\mathrm{Fr}(\Sigma(M))$ are essential in $M$, the region of parallelism cannot contain any solid torus components of $\Sigma(M)$. Therefore $F$ is also inessential in $\widehat{M}$, which is a contradiction since it is a component of $\mathrm{Fr}(\Sigma(\widehat{M}))$. If $F$ is an inessential annulus in $M$, then it together with an annulus in $\partial M$ bounds a solid torus in $M$, which cannot contain any solid torus components of $\Sigma(M)$. Therefore $F$ is inessential in $\widehat{M}$, which again is contradictory. This completes the proof of the claim.

It follows that each component $R$ of the closure of $\Sigma(\widehat{M})-(\Sigma(M)-\mathcal{N}(\Delta))$ is an essential fibered submanifold of $M$. Therefore $R$ is admissibly isotopic in $M$ into $\Sigma(M)$. By Lemma 6.2, $R$ is a product region whose ends are equal to the components of its frontier. We claim that exactly one of the components of $\mathrm{Fr}(R)$ lies in $\Sigma(M)$.

Suppose neither does, so that $R$ is a component of $\Sigma(\widehat{M})$ which does not meet $\Sigma(M)-\mathcal{N}(\Delta)$. Since $R$ is isotopic in $M$ into $\Sigma(M)$, each component of the frontier of $R$ is parallel in $\overline{M-\Sigma(M)}$, and hence in $\widehat{M}$, to a component of $\mathrm{Fr}(\Sigma(M))$. Since $R$ is a product, it follows that $R$ is admissibly homotopic in $\widehat{M}$ into $\Sigma(\widehat{M})-R$. Therefore, $\Sigma(\widehat{M})-R$ has the engulfing property in $\widehat{M}$, contradicting the minimality of $\Sigma(\widehat{M})$.

Suppose now that both components of $\mathrm{Fr}(R)$ lie in $\Sigma(M)$. Since $\widehat{M}$ is not an $I$-bundle over the torus, the fibering of $\Sigma(\widehat{M})$ must be a Seifert fibering on the components that meet $\partial \overline{\mathcal{N}(\Delta)}$, so it extends over $\overline{\mathcal{N}(\Delta)}$. Since $\mathrm{Fr}(R\cup\Sigma(M))$ consists of essential annuli, and $R\cup\Sigma(M)$ is contained in the topological interior of $\Sigma(\widehat{M})\cup\mathcal{N}(\Delta)$, Corollary 5.7 of Johannson [19] shows it is an essential fibered submanifold of $M$. Since $R\cup\Sigma(M)$ also has the engulfing property, this contradicts the uniqueness of $\Sigma(M)$.

We conclude that each component of the closure of $\Sigma(\widehat{M})-(\Sigma(M)-\mathcal{N}(\Delta))$ is a product region with one end a component of $\mathrm{Fr}(\Sigma(M))$ and the other a component of $\mathrm{Fr}(\Sigma(\widehat{M}))$. It follows that $\Sigma(M)-\mathcal{N}(\Delta)$ is isotopic in $\widehat{M}$ to $\Sigma(\widehat{M})$, so is characteristic.

Since $\Sigma(M)-\mathcal{N}(\Delta)$ is characteristic in $\widehat{M}$, every incompressible torus in $\widehat{M}$ is homotopic in $\widehat{M}$ into $\Sigma(M)-\mathcal{N}(\Delta)$, hence is homotopic into a torus boundary component of $\widehat{M}$. This shows that $\widehat{M}$ is atoroidal and hence hyperbolizable. □

Suppose now that $\widehat{M}$ is constructed as in Lemma 10.2, and that all the solid tori containing the loops $\Delta$ are primitive. Choose a meridian-longitude system $(m_i, l_i)$ for the boundary $T_i$ of each component $\mathcal{N}(\delta_i)$ of $\mathcal{N}(\Delta)$ so that the longitude is parallel in $M - \mathcal{N}(\delta_i)$ into $\mathcal{V} \cap \partial M$, and so that the meridian bounds a disc in $\overline{\mathcal{N}(\delta_i)}$. The resulting meridian-longitude system is said to be *natural*.

In what follows, we often consider the specific Dehn surgery coefficients $(\mathbf{p}_n, \mathbf{q}_n) = ((1, n), \ldots, (1, n))$ and the resulting manifold $\widehat{M}(\mathbf{p}_n, \mathbf{q}_n)$. For this case, Lemma 10.3 below assures us that $\widehat{M}(\mathbf{p}_n, \mathbf{q}_n)$ is homeomorphic to $M$.

We also wish to exhibit an explicit compact core for $\widehat{M}(\mathbf{p}_n, \mathbf{q}_n)$. Let $V_i$ be the element of $\mathcal{V}$ containing $\mathcal{N}(\delta_i)$ and let $A_i$ be an essential annulus in $V_i - \mathcal{N}(\delta_i)$ with one boundary component parallel to the longitude $l_i$ on $T_i$ and the other in $V_i \cap \partial M$. We form $M_0$ by removing an open regular neighborhood of $A_i$ from $\widehat{M}$. Then, since each component $W_i$ of $M - M_0$ is a solid torus intersecting $M_0$ in an annulus which is primitive in $W_i$, there exists a strong deformation retraction $\tau \colon M \to M_0$. We call a submanifold of $M$ constructed in this manner, a *standard* compact core for $M$.

LEMMA 10.3. *Let $M$ be a compact, hyperbolizable 3-manifold with non-empty incompressible boundary, let $\mathcal{V}$ be a collection of primitive solid torus components of the characteristic submanifold $\Sigma(M)$ of $M$, let $\Delta$ be the collection of core curves of tori in $\mathcal{V}$, and let $\mathcal{N}(\Delta)$ be an open regular neighborhood of $\Delta$. Let $\{T_1, \ldots, T_k\}$ be the boundary components of $\mathcal{N}(\Delta)$ and let $(m_i, l_i)$ be natural meridian-longitude systems. Then, for any $n$ and for $(\mathbf{p}_n, \mathbf{q}_n) = ((1, n), \ldots, (1, n))$, we have that $\widehat{M}(\mathbf{p}_n, \mathbf{q}_n)$ is homeomorphic to $M$. Moreover, if $M_0$ is a standard compact core for $M$, then $i_n(M_0)$ is a compact core for $\widehat{M}(\mathbf{p}_n, \mathbf{q}_n)$ where $i_n$ denotes the inclusion of $\widehat{M}$ into $\widehat{M}(\mathbf{p}_n, \mathbf{q}_n)$.*

*Proof of Lemma* 10.3. Let $A_i$ be the annuli used in the construction of the standard compact core $M_0$. Let $h_n \colon \widehat{M} \to \widehat{M}$ be obtained as a composition of $n$-fold Dehn twists about each $A_i$. We may assume that $h_n$ agrees with the identity map on $M_0$. Then $h_n$ takes a $(1, 0)$-curve in $T_i$ to a $(1, n)$-curve in $T_i$. Thus, $h_n$ extends to a homeomorphism $H_n \colon \widehat{M}(\mathbf{p}_0, \mathbf{q}_0) \to \widehat{M}(\mathbf{p}_n, \mathbf{q}_n)$. Since $\widehat{M}(\mathbf{p}_0, \mathbf{q}_0) = M$, $H_n$ is the desired homeomorphism. Because $M_0$ is a compact core for $M$, $H_n(M_0)$ is a compact core for $\widehat{M}(\mathbf{p}_n, \mathbf{q}_n)$. Moreover, since $h_n$ agrees with the identity map on $M_0$, $H_n(M_0) = i_n(M_0)$. □

## 11. Deformation spaces which go bump in the night

In this section, we prove Theorem B. It generalizes the main theorem of [4], and its proof follows much the same outline. The key new ingredient is Proposition 9.1. It assures that the construction in [4], which deals only with a very special case, can be carried out in our much more general situation.

THEOREM B. *Let $M$ be a compact, hyperbolizable 3-manifold with non-empty incompressible boundary, and let $[(M_1, h_1)]$ and $[(M_2, h_2)]$ be two elements of $\mathcal{A}(M)$. If $[(M_2, h_2)]$ is obtained from $[(M_1, h_1)]$ by applying a primitive shuffle, then the associated components of $MP(\pi_1(M))$ have intersecting closures.*

*Proof of Theorem* B. To set notation, let $U_j$ be the component of $MP(\pi_1(M))$ corresponding to $[(M_j, h_j)] \in \mathcal{A}(M)$, so that $\Psi^{-1}([(M_j, h_j)]) \cap MP(\pi_1(M)) = U_j$. We need to show that there exists a representation $\rho$ in $AH(\pi_1(M))$ in the intersection $\overline{U_1} \cap \overline{U_2}$ of the closures of the two components.

Let $\phi\colon M_1 \to M_2$ be a primitive shuffle with respect to $\mathcal{V}_1$ and $\mathcal{V}_2$ such that $\phi \circ h_1$ is homotopic to $h_2$. Let $\Delta = \{\delta_1, \ldots, \delta_k\}$ be the collection of core curves of the solid tori $\{V_2^1, \ldots, V_2^k\}$ of $\mathcal{V}_2$, let $\mathcal{N}(\Delta)$ be an open regular neighborhood of $\Delta$ in $\mathcal{V}_2$ chosen so that the closure of $\mathcal{N}(\Delta)$ lies in the interior of $\mathcal{V}_2$, and let $\widehat{M}_2 = M_2 - \mathcal{N}(\Delta)$. By Proposition 9.1, we may assume that $\phi(M_1) \subset \widehat{M}_2$ and that $\phi$ lifts to an embedding in some cover of $\widehat{M}_2$. Notice that $\widehat{M}_2$ is compact by construction and is hyperbolizable by Lemma 10.2. Let $T = \{T_1, \ldots, T_k\}$ be the collection of torus boundary components of $\widehat{M}_2$ coming from the elements of $\mathcal{N}(\Delta)$, so that $T_i$ is parallel to $\partial V_2^i$ in $V_2^i$.

Choose a meridian-longitude system $(m_i, l_i)$ for each $T_i$ so that the longitude is parallel, in $\widehat{M}_2$, to a core curve of a component of the frontier of $V_2^i$ in $M_2$ and so that the meridian bounds a disc in $\mathcal{N}(\delta_i)$; such a system was referred to as a *natural meridian-longitude system* in Section 10. For each $1 \leq i \leq k$, let $A_i$ be an essential annulus in $W_i = V_2^i - \mathcal{N}(\delta_i)$ with one boundary component the longitude $l_i$ on $T_i$ and the other boundary component lying in $\partial V_2^i \cap \partial M_2$. Let $M_2^0$ be obtained from $\widehat{M}_2$ by removing, for all $1 \leq i \leq k$, an open regular neighborhood of $A_i$ which is itself contained in $W_i$. We observed in Section 10 that $M_2^0$ is a compact core for $M_2$ and that there exists a strong deformation retraction $\tau\colon M_2 \to M_2^0$.

Since $\widehat{M}_2$ is compact and hyperbolizable, Thurston's geometrization theorem (see [33]) implies that there exists a geometrically finite hyperbolic 3-manifold $\widehat{N} = \mathbf{H}^3/\widehat{\Gamma}$ such that every parabolic element of $\widehat{\Gamma}$ lies in a rank two parabolic subgroup of $\widehat{\Gamma}$, and an orientation-preserving homeomorphism $\psi\colon \operatorname{int}(\widehat{M}_2) \to \widehat{N}$. Let $\{\beta_n\colon \widehat{\Gamma} \to \operatorname{PSL}_2(\mathbf{C})\}$ be the sequence of representations given by applying Theorem 10.1 to $\widehat{\Gamma}$ with the sequence of $k$-tuples $\{(\mathbf{p}_n, \mathbf{q}_n) = ((1, n), \ldots, (1, n))\}$. Let $\widehat{M}_2(n) = \widehat{M}_2(\mathbf{p}_n, \mathbf{q}_n)$, let $i_n\colon \widehat{M}_2 \to \widehat{M}_2(n)$ be the inclusion map, and let $\psi_n\colon \operatorname{int}(\widehat{M}_2(n)) \to N_n = \mathbf{H}^3/\beta_n(\widehat{\Gamma})$ be the homeomorphism from the interior of $\widehat{M}_2(n)$ to $N_n$ with the property that $(\psi_n)_* \circ (i_n)_*$ is conjugate to $\beta_n \circ \psi_*$. By Lemma 10.3, $\widehat{M}_2(n)$ is homeomorphic to $M_2$ and $i_n(M_2^0)$ is a compact core for $\widehat{M}_2(n)$.

We now observe that $i_n \circ \phi$ is itself a primitive shuffle equivalence. We begin by noting that the composition $i_n \circ \phi$ makes sense, since $\phi(M_1) \subset \widehat{M}_2 \subset M_2$, and that $i_n \circ \phi$ restricts to an orientation-preserving embedding of $\overline{M_1 - \mathcal{V}_1}$ into $\widehat{M}_2(n)$. Moreover,
$$\mathcal{V}_2(n) = \widehat{M}_2(n) - i_n(\phi(M_1 - \mathcal{V}_1))$$
is a collection of solid tori, since $\mathcal{V}_2(n)$ is obtained from $\mathcal{V}_2$ by $(1, n)$-Dehn filling along the core curve of each solid torus component of $\mathcal{V}_2$. Each component of the frontier of $\mathcal{V}_2(n)$ in $\widehat{M}_2(n)$ is the image under $i_n$ of a primitive essential annulus in $M_2^0$. Since $i_n(M_2^0)$ is a compact core for $\widehat{M}_2(n)$, we may conclude that each component of the frontier of $\mathcal{V}_2(n)$ is a primitive essential annulus in $M_2(n)$, and hence that $\mathcal{V}_2(n)$ is a collection of primitive essentially embedded solid tori. Since $i_n \circ \phi$ is a homeomorphism from the frontier of $\mathcal{V}_1$ to the frontier of $\mathcal{V}_2(n)$, $i_n \circ \phi$ is a homotopy equivalence from $\mathcal{V}_1$ to $\mathcal{V}_2(n)$. Lemma 5.2 implies that $i_n \circ \phi$ is a shuffle homotopy equivalence with respect to $\mathcal{V}_1$ and $\mathcal{V}_2(n)$ for all $n$. Proposition 6.1 assures us that we may assume that each element of $\mathcal{V}_2(n)$ is a solid torus component of $\Sigma(\widehat{M}_2(n))$, so that $i_n \circ \phi$ is indeed a primitive shuffle equivalence.

One may similarly show that $i_n \circ \tau \circ \phi$ is a primitive shuffle, with respect to $\mathcal{V}_1$ and $\mathcal{V}_2(n)$, from $M_1$ to $\widehat{M}_2(n)$. Since $i_n \circ \tau \circ \phi$ and $i_n \circ \phi$ agree on $M_1 - \mathcal{V}_1$, Lemma 7.1 implies that there exists an orientation-preserving homeomorphism $r_n \colon M_1 \to M_1$, which is the identity on $M_1 - \mathcal{V}_1$, such that $i_n \circ \tau \circ \phi \circ r_n$ is homotopic to $i_n \circ \phi$.

For each $n$, define $\rho'_n = \beta_n \circ \psi_* \circ \phi_*$. Since $\beta_n \circ \psi_*$ is conjugate to $(\psi_n)_* \circ (i_n)_*$, we see that $\rho'_n$ is conjugate to $(\psi_n)_* \circ (i_n)_* \circ \phi_*$. Since $i_n \circ \phi$ is a homotopy equivalence and $\psi_n$ is a homeomorphism, $\rho'_n$ is a discrete faithful representation with image $\pi_1(N_n) = \beta_n(\widehat{\Gamma})$. The first part of Theorem 10.1 then implies that $\rho'_n \in MP(\pi_1(M_1))$. Moreover, $\{\rho'_n\}$ converges to $\rho' = \psi_* \circ \phi_*$.

As $M_1$ is a compact, hyperbolizable 3-manifold with nonempty incompressible boundary, we may consider the map $\Psi_1 \colon AH(\pi_1(M_1)) \to \mathcal{A}(M_1)$. Let $i'_n$ be an embedding of $\widehat{M}_2$ into the interior of $\widehat{M}_2(n)$ which is isotopic to $i_n$. Then $\psi_n(i'_n(M_2^0))$ is a compact core of $N_n$. Moreover, $s_n = \psi_n \circ i'_n \circ \tau \circ \phi \circ r_n$ is a homotopy equivalence from $M_1$ to $N_n$ with image in $\psi_n(i'_n(M_2^0))$ such that $(s_n)_*$ is conjugate to $\rho'_n$. Therefore,
$$\Psi_1(\rho'_n) = [(\psi_n(i'_n(M_2^0)), s_n)] = [(M_2^0, \tau \circ \phi \circ r_n)] = [(M_2, \phi \circ r_n)]$$
for all $n$; the last equality follows from the observation that $\tau$ is homotopic to an orientation-preserving homeomorphism from $M_2$ to $M_2^0$. However, since $\phi \circ r_n$ and $\phi$ are both primitive shuffle equivalences with respect to $\mathcal{V}_1$ and $\mathcal{V}_2$ which agree on $M_1 - \mathcal{V}_1$, Lemma 7.1 implies that, for all $n$, there exists an orientation-preserving homeomorphism $r'_n \colon M_2 \to M_2$ such that $r'_n \circ \phi$ is homotopic to $\phi \circ r_n$. Therefore,
$$\Psi_1(\rho'_n) = [(M_2, \phi \circ r_n)] = [(M_2, r'_n \circ \phi)] = [(M_2, \phi)]$$
for all $n$.

It is easy to construct a map $\overline{\phi}\colon M_1 \to M_2$, homotopic to $\phi$, so that $\overline{\phi}(M_1) \subset \operatorname{int}(\widehat{M_2})$ and $\overline{\phi}$ still lifts to an embedding of $M_1$ into some cover of $\operatorname{int}(\widehat{M_2})$. Then $\psi \circ \overline{\phi}$ lifts to an orientation-preserving embedding $g\colon M_1 \to N$ of $M_1$ into the cover $N = \mathbf{H}^3/\rho'(\pi_1(M_1))$ of $\widehat{N}$ associated to $\psi_*(\phi_*(\pi_1(M_1))) = \rho'(\pi_1(M_1))$. Since $g_* = \rho'$,
$$\Psi_1(\rho') = [(g(M_1), g)] = [(M_1, \operatorname{id})].$$

Let $\rho = \rho' \circ (h_1)_*$ and let $\rho_n = \rho'_n \circ (h_1)_*$. Then, for all $n$, we have that $\rho_n \in MP(\pi_1(M))$ and that $\Psi(\rho_n) = [(M_2, \phi \circ h_1)] = [(M_2, h_2)]$. In particular, since $\rho = \lim \rho_n$, we see that $\rho$ lies in the closure of the component $U_2$ of $MP(\pi_1(M))$ corresponding to $[(M_2, h_2)]$. Moreover, $\Psi(\rho) = [(M_1, h_1)]$.

It remains to show that $\rho$ lies in the closure of the component $U_1$ of $MP(\pi_1(M))$ corresponding to $[(M_1, h_1)]$. Since $\rho'(\pi_1(M_1))$ is a finitely generated subgroup of a co-infinite volume, geometrically finite Kleinian group, it is itself geometrically finite (see Proposition 7.1 in [33]). Corollary 6 of Ohshika [35] then guarantees that we may write $\rho' = \lim \alpha'_n$, where $\alpha'_n \in MP(\pi_1(M_1))$ and $\Psi_1(\alpha'_n) = \Psi_1(\rho') = [(M_1, \operatorname{id})]$ for all $n$. Let $\alpha_n = \alpha'_n \circ (h_1)_*$. Then, for all $n$, we have that $\alpha_n \in MP(\pi_1(M))$ and $\Psi(\alpha_n) = [(M_1, h_1)]$. Since $\rho = \lim \alpha_n$, we see that $\rho$ lies in the closure of $U_1$.

Hence, we see that $\rho$ lies in the intersection $\overline{U_1} \cap \overline{U_2}$ of the closures of the two components of $MP(\pi_1(M))$ associated to $[(M_1, h_1)]$ and $[(M_2, h_2)]$, as desired. $\square$

*Remark.* One may also use Theorem 10.1 to construct the sequence $\alpha'_n$ above. Let $\widehat{M_1}$ be obtained from $M_1$ by removing regular neighborhoods of all the core curves of elements of $\mathcal{V}_1$. Let $M_1^0$ be a standard compact core for $M_1$ constructed as above and let $\tau_1\colon M_1 \to M_1^0$ be a strong deformation retraction. One may use the Klein-Maskit combination theorems to construct a Kleinian group $\widehat{\Gamma}_1$ uniformizing $\widehat{M_1}$ and containing $\rho(\pi_1(M))$ as a precisely embedded subgroup. Moreover, we may construct $\widehat{\Gamma}_1$ so that there exists a homeomorphism $\psi_1\colon \operatorname{int}(\widehat{M_1}) \to \mathbf{H}^3/\widehat{\Gamma}$ such that $\rho = (\psi_1)_* \circ (\tau_1)_*$. Let $\{\widehat{\alpha}_n\colon \widehat{\Gamma}_1 \to \operatorname{PSL}_2(\mathbf{C})\}$ be the sequence of representations given by applying Theorem 10.1 to $\widehat{\Gamma}_1$ with the sequence of $k$-tuples $\{(1, n), \ldots, (1, n)\}$. Then $\{\alpha'_n = \widehat{\alpha}_n \circ (\psi_1)_* \circ (\tau_1)_*\}$ is the desired sequence.

## 12. Proofs of corollaries

In this section we give proofs of the corollaries stated in the introduction.

COROLLARY 1. *Let $M$ be a compact, hyperbolizable $3$-manifold with nonempty incompressible boundary. If $U_1$ and $U_2$ are components of $MP(\pi_1(M))$, then $\overline{U}_1 \cap \overline{U}_2$ is nonempty if and only if $\Psi(U_1)$ and $\Psi(U_2)$ are primitive shuffle equivalent.*

*Proof of Corollary* 1. If $\Psi(U_1)$ and $\Psi(U_2)$ are primitive shuffle equivalent, then Theorem B implies that $\overline{U}_1 \cap \overline{U}_2$ is nonempty. On the other hand, if $\rho \in \overline{U}_1 \cap \overline{U}_2$, then Theorem A implies that $\Psi(\rho)$ is primitive shuffle equivalent to both $\Psi(U_1)$ and $\Psi(U_2)$. Therefore, since primitive shuffle equivalence is an equivalence relation, $\Psi(U_1)$ and $\Psi(U_2)$ are primitive shuffle equivalent. □

COROLLARY 2. *Let $M$ be a compact, hyperbolizable 3-manifold with nonempty incompressible boundary. Then, the components of $\mathrm{MP}(\pi_1(M))$ cannot accumulate in $AH(\pi_1(M))$. In particular, the closure $\overline{\mathrm{MP}(\pi_1(M))}$ of $\mathrm{MP}(\pi_1(M))$ is the union of the closures of the components of $\mathrm{MP}(\pi_1(M))$.*

*Proof of Corollary* 2. Suppose that the components of $\mathrm{MP}(\pi_1(M))$ accumulated in $AH(\pi_1(M))$. Then there would exist a sequence $\{\rho_n\}$ in $\mathrm{MP}(\pi_1(M))$, converging to some $\rho$ in $AH(\pi_1(M))$, with $\Psi(\rho_n) \neq \Psi(\rho_m)$ for $n \neq m$. By Theorem A, $\Psi(\rho_n)$ is primitive shuffle equivalent to $\Psi(\rho)$ for all sufficiently large $n$. Thus, $\Psi(\rho_n) \in q^{-1}(\widehat{\Psi}(\rho))$ for all sufficiently large $n$. However, Proposition 7.2 implies that $q^{-1}(\widehat{\Psi}(\rho))$ is finite. This contradiction establishes our first claim.

Since for any locally finite collection $\{A_\alpha\}$ of subsets of a space $X$, the closure of $\cup A_\alpha$ is the union of the closures of the $A_\alpha$, the second statement then follows. □

COROLLARY 3. *Let $M$ be a compact, hyperbolizable 3-manifold with nonempty incompressible boundary. Then, the components of $\overline{\mathrm{MP}(\pi_1(M))}$ are in one-to-one correspondence with the elements of $\widehat{\mathcal{A}}(M)$.*

*Proof of Corollary* 3. Let $y \in \widehat{\mathcal{A}}(M)$. Since the restriction of $\widehat{\Psi}$ to $\mathrm{MP}(\pi_1(M))$ is surjective, there exists a component $U$ of $\mathrm{MP}(\pi_1(M))$ with $U \subseteq \widehat{\Psi}^{-1}(y)$. Let $C$ be the connected component of $\overline{\mathrm{MP}(\pi_1(M))}$ that contains $U$. Since, by Corollary 8.1, $\widehat{\Psi}$ is locally constant, $C \subseteq \widehat{\Psi}^{-1}(y)$. Now suppose that $\rho \in \overline{\mathrm{MP}(\pi_1(M))}$ and $\widehat{\Psi}(\rho) = y$. By Corollary 2, there is a component $V$ of $\mathrm{MP}(\pi_1(M))$ with $\rho \in \overline{V}$. Again since $\widehat{\Psi}$ is locally constant, $\widehat{\Psi}(\overline{V}) = y$. By Theorem B, $\overline{V} \cap \overline{U}$ is nonempty, so that $\overline{V} \subseteq C$ and hence $\rho \in C$. Therefore $\widehat{\Psi}^{-1}(y) = C$. Since $\widehat{\Psi}$ is surjective, this exhibits an explicit one-to-one correspondence between components of $\overline{\mathrm{MP}(\pi_1(M))}$ and elements of $\widehat{\mathcal{A}}(M)$. □

COROLLARY 4. *Let $M$ be a compact, hyperbolizable 3-manifold with nonempty incompressible boundary. Then, $\overline{\mathrm{MP}(\pi_1(M))}$ has infinitely many components if and only if $M$ has double trouble. Moreover, if $M$ has double trouble, then $AH(\pi_1(M))$ has infinitely many components.*

*Proof of Corollary* 4. It follows immediately from Corollary 3 that $\overline{\mathrm{MP}(\pi_1(M))}$ has infinitely many components if and only if $\widehat{\mathcal{A}}(M)$ has infinitely many elements. The fact that $q: \mathcal{A}(M) \to \widehat{\mathcal{A}}(M)$ is finite-to-one implies that $\widehat{\mathcal{A}}(M)$ has infinitely many elements if and only if $\mathcal{A}(M)$ has infinitely many elements. The results of [12] imply that $\mathcal{A}(M)$ has infinitely many elements if and only if $M$ has double trouble. This establishes our first claim.

Now suppose that $M$ has double trouble, so that $\widehat{\mathcal{A}}(M)$ has infinitely many elements. Corollary 8.1 gives that $\widehat{\Psi}\colon AH(\pi_1(M)) \to \widehat{\mathcal{A}}(M)$ is locally constant. Since $\widehat{\Psi}$ is also surjective, $AH(\pi_1(M))$ has infinitely many components as claimed. □

COROLLARY 5. *Let $M$ be a compact, hyperbolizable 3-manifold with nonempty incompressible boundary. If $q\colon \mathcal{A}(M) \to \widehat{\mathcal{A}}(M)$ is not injective, then $\overline{MP(\pi_1(M))}$ and $AH(\pi_1(M))$ are not manifolds.*

*Proof of Corollary* 5. Suppose that $\overline{MP(\pi_1(M))}$ is a manifold and $q$ is not injective. Theorem B implies that there exist distinct components $U$ and $V$ of $MP(\pi_1(M))$ and a point $\rho$ in the closure of both $U$ and $V$. Let $Y$ denote the component of $\overline{MP(\pi_1(M))}$ which contains both $U$ and $V$. Let $Z$ denote the set of representations $\rho$ in $\overline{MP(\pi_1(M))}$ with the property that there exists a core curve $\alpha$ of a primitive solid torus component of $\Sigma(M)$ with $\rho(\alpha)$ parabolic. Corollary 8.2 implies that $\Psi$ is continuous on $AH(\pi_1(M)) - Z$, so that $Z$ must disconnect $Y$. Hence, if $n = \dim(MP(\pi_1(M))) = \dim(X_T(\pi_1(M))) = \dim(Y)$, then $Z$ must have topological dimension at least $n - 1$ (see Theorem IV.4 in Hurewicz-Wallman [17]). On the other hand, since $Z$ lies in a finite union of (complex) codimension one subvarieties of $X_T(\pi_1(M))$, $Z$ must have topological dimension at most $n - 2$. This contradiction establishes that $\overline{MP(\pi_1(M))}$ is not a manifold if $q$ is not injective.

The argument which proves that $AH(\pi_1(M))$ is not a manifold is exactly the same. □

UNIVERSITY OF SOUTHAMPTON, SOUTHAMPTON, ENGLAND
*E-mail address*: jwa@maths.soton.ac.uk

UNIVERSITY OF MICHIGAN, ANN ARBOR, MI
*E-mail address*: canary@math.lsa.umich.edu

UNIVERSITY OF OKLAHOMA, NORMAN, OK
*E-mail address*: dmccullough@math.ou.edu